\documentclass[10pt]{article}
\usepackage[dvips]{graphicx}
\usepackage{amsmath,amsfonts,amssymb,latexsym,epsfig}
\include{references}
\usepackage{mathrsfs}
\usepackage{verbatim}
\usepackage{latexsym}
\usepackage{amsthm}
\usepackage{amssymb}
\usepackage{graphics}
\usepackage{amsbsy}
\usepackage{enumerate}
\usepackage{times}
\newtheorem{theorem}{Theorem}[section]
\newtheorem{lemma}[theorem]{Lemma}
\newtheorem{remark}[theorem]{Remark}
\newtheorem{corollary}[theorem]{Corollary}
\newtheorem{prop}[theorem]{Proposition}
\newtheorem{assumptions}[theorem]{Assumptions}

\numberwithin{equation}{section}
\graphicspath{{figures/}}
\newcommand{\E}{{\mathbb E}}
\newcommand{\bbE}{{\mathbb E}}

\newcommand{\LL}{{\mathcal L}}

\newcommand{\T}{{\mathbb T}}
\newcommand{\R}{{\mathbb R  }}

\newcommand{\bbT}{\mathbb{T}}
\newcommand{\bbR}{\mathbb{R}}

\newcommand{\eps}{\epsilon}

\newcommand{\la}{\langle}
\newcommand{\ra}{\rangle}

\begin{document}
\setlength{\baselineskip}{10pt}
\title{PARAMETER ESTIMATION FOR MULTISCALE DIFFUSIONS}
\author{G.A. Pavliotis\footnote{Corresponding author. 
E-mail address: g.paviotis@maths.warwick.ac.uk.} \\
        Department of Mathematics\\
    Imperial College London \\
        London SW7 2AZ, UK \\
        and \\
        A.M. Stuart\footnote{E-mail address: stuart@maths.warrwick.ac.uk.} \\
        Mathematics Institute \\
        Warwick University \\
        Coventry CV4 7AL, UK
                    }
\maketitle

\begin{abstract}
We study the problem of parameter estimation for time-series possessing two, widely
separated, characteristic time scales. The aim is to understand situations where it is
desirable to fit a homogenized singlescale model to such multiscale data. We
demonstrate, numerically and analytically, that if the data is sampled too finely then
the parameter fit will fail, in that the correct parameters in the homogenized model are
not identified. We also show, numerically and analytically, that if the data is
subsampled at an appropriate rate then it is possible to estimate the coefficients of the
homogenized model correctly.
\end{abstract}


%
%%%%%%%%%%%%%%%%%%%%%%%%%%%%%%%%%%%%%%%%%%%%%%%%%%%%%%%%%%%%%%%%%%%%%%%%%%%%%%%%%%%%%%%%%
%
%                                      INTRODUCTION
%
%%%%%%%%%%%%%%%%%%%%%%%%%%%%%%%%%%%%%%%%%%%%%%%%%%%%%%%%%%%%%%%%%%%%%%%%%%%%%%%%%%%%%%%%%%
%
\section{Introduction}
\label{sec:intro}

Parameter estimation for continuous time stochastic models is an increasingly important
part of the overall modelling strategy in a wide variety of applications. It is quite
often the case that the data to be fitted to a diffusion process has a multiscale
character. One example is the field of  molecular dynamics, where it is desirable to find
effective models for low dimensional phenomena (such as conformational dynamics, vacancy
diffusion and so forth) which are embedded within higher dimensional time-series. Another
example is the ocean--atmosphere sciences where it is desirable to find effective models
for large--scale structures, whilst representing the small--scales stochastically.  The
multiscale structure of the data in these problems renders the problem of parameter
estimation very subtle, and great care has to be taken in order to estimate the
coefficients correctly. The aim of the paper is to shed light on this estimation problem
through the study of a simple class of model problems, typical of those arising in
molecular dynamics.

In econometrics and finance, the problem of estimating parameters for continuous time
diffusion processes in the presence of small scale fluctuations (market microstructure
noise) has been considered by A\"{i}t--Sahalia and collaborators
\cite{AitMykZha05b,AitMykZha05a} and more recently in \cite{BaNiHaLuSh06}. In that work
the microscale is input as an independent white observational noise that is superimposed
on--top of a singlescale diffusion process. We have a somewhat different framework: we
work in the context of coupled systems of diffusions exhibiting multiple scales. Our aim
is to fit a singlescale homogenized diffusion to data. Models similar to the ones
considered in this paper have been studied extensively in finance, see \cite{fouque00}
and the reference therein. In that book there is discussion of parameter estimation for
multiscale diffusions, with emphasis on the estimation of the rate of mean reversion of
volatility from historical asset price data; see \cite[Ch. 4]{fouque00}.

Various numerical algorithms for
diffusions with multiple scales have been developed \cite{Vand03} and analyzed
\cite{ELV05}. Those papers are finely honed to optimize the fitting of the homogenized
diffusion in situations where the multiscale model is known explicitly. In contrast, in
this paper we introduce multiscale diffusions primarily as a device to generate
multiscale data; we do not assume that the multiscale model is available to us when doing
parameter estimation. This enables us to gain understanding of parameter estimation in
situations where the multiscale data is given to us from experiments, or comes from a
model where the scale--separation is not explicit. Two recent papers contain numerical
experiments relating to the extraction of averaged or homogenized diffusions from data
generated by a multiscale diffusion; see \cite{Cald06,CromVanEij06b}.

Despite differences from the framework used in
\cite{AitMykZha05b,AitMykZha05a,BaNiHaLuSh06} to study problems arising
in econometrics and finance, similarities with our work remain:
trying to fit the
models on the basis of data sampled at too high a frequency leads to incorrect parameter
inference; furthermore, there is an optimal subsampling rate for the data to obtain
correct inference.

There are two forms of multiscale diffusions which are of particular
interest in the context of parameter estimation. The first gives rise
to {\bf averaging} for SDEs, and the second to {\bf homogenization}
for SDEs. For averaging one has, for $\eps \ll 1$,
\begin{subequations}
\begin{eqnarray}
 d x^\eps(t) &=& f(x^\eps(t),y^\eps(t)) \, dt+\alpha(x^\eps(t),y^\eps(t)) \,  dU(t),  \\
 d y^\eps(t) &=& \frac{1}{\eps} g(x^\eps(t),y^\eps(t)) \, dt+\frac{1}{\sqrt\eps}
 \beta(x^\eps(t),y^\eps(t)) \, dV(t),
\end{eqnarray}
\label{e:averg}
\end{subequations}
with $U,V$ standard Brownian motions. Averaging $f$ and $\alpha \alpha^T$
over the invariant measure of the $y^\eps$ equation, with $x^\eps$ viewed as fixed,
gives an averaged SDE for $x$. The fast process $y$, with timescale $\eps$,
is eliminated. For homogenization one has
\begin{subequations}
\begin{eqnarray}
d x^\eps(t) &=&  \left( \frac{1}{\eps} f_0(x^\eps(t),y^\eps(t)) +
f_1(x^\eps(t),y^\eps(t)) \right) dt \nonumber \\
&+& \alpha(x^\eps(t),y^\eps(t)) \, dU( t),\\
d y^\eps(t) &=& \frac{1}{\eps^2} g(x^\eps(t),y^\eps(t)) \, dt +
\frac{1}{\eps} \beta(x^\eps(t),y^\eps(t)) \, dV(t),
\end{eqnarray}
\label{e:homog}
\end{subequations}
where it is assumed that $f_0$ averages to zero against the invariant measure
of the fast process $y^\eps$ with $x^\eps$ fixed.
Now $y^\eps$ has time-scale $\eps^2$ and is eliminated.
The fluctuations in $f_0$, suitably amplified by $\eps^{-1}$, induce ${\cal O}(1)$ effects
in the homogenized equation for $x^\eps$. In both cases \eqref{e:averg} and \eqref{e:homog} it is
possible to show \cite{lions} that the process $x^\eps(t)$ converges in law, as $\eps
\rightarrow 0$, to the solution of an effective SDE of the form
\begin{equation}\label{e:effect}
d x(t) = F(x(t)) dt + A(x(t)) d U(t).
\end{equation}
Explicit formulae can be derived for the effective coefficients $F(x)$ and $A(x)$ in the
above equation \cite{lions, PavlSt06b}. A natural question that arises then is how to fit
an SDE of the form \eqref{e:effect} to data generated by a multiscale stochastic equation
of the form \eqref{e:averg} or \eqref{e:homog}, under the assumption of scale separation,
i.e. when $\eps \ll 1$. This paper is a first attempt towards the study of this
interesting problem, for a specific class of SDEs of the form \eqref{e:homog}.

Our basic model will be the first order Langevin equation
\begin{equation}
d x^\eps(t) = - \nabla V \left(x^\eps(t), \frac{x^\eps(t)}{\eps}; \alpha \right)  dt +
\sqrt{2 \sigma } d \beta(t),
\label{e:main}
\end{equation}
where $\beta(t)$ denotes standard Brownian motion on $\R^d$ and $\sigma$ is a positive
constant. The two--scale potential $ V^\eps \left(x, y; \alpha \right)$ is assumed to
consist of a large--scale and a fluctuating part
\begin{equation}
V ( x, y ; \alpha) = \alpha V(x) + p(y).
\label{e:potential}
\end{equation}
As we show explicitly in \eqref{e:eqns_motion} this set-up puts us in
the framework of homogenization for SDEs.

\begin{figure}
\begin{center}
\includegraphics[width=3.0in, height = 3.0in]{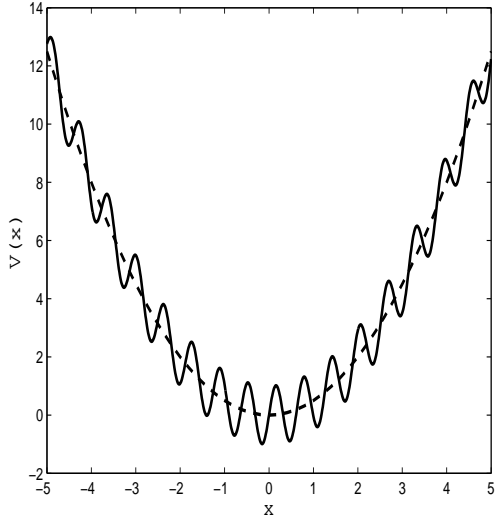}
\caption{$V^\eps(x, x/\eps) = \frac{1}{2} x^2 + \sin \left( \frac{x}{\eps} \right)$
with $ \eps = 0.1$ and averaged potential $V(x) = \frac{1}{2} x^2$.} \label{fig:potential}
\end{center}
\end{figure}

Under \eqref{e:potential}, the SDE \eqref{e:main} becomes
\begin{equation}
d x^\eps(t) = - \alpha \nabla V(x^\eps(t)) \, dt - \frac{1}{\eps}\nabla p \left(
\frac{x^\eps(t)} {\eps} \right) \, dt + \sqrt{2 \sigma} \, d \beta (t).
\label{e:xeps_V}
\end{equation}
If $p$ is periodic on $\bbT^d$ and sufficiently smooth, then it is well
known (see \cite{lions, pardoux} for example) that, as $\eps \rightarrow 0$,
the solution $x^\eps(t)$ of $\eqref{e:main}$ converges in law to the
solution of the SDE
\begin{equation}
 d x(t) = -\alpha K \nabla V(x(t)) dt + \sqrt{2 \sigma K} d \beta (t),
\label{e:lim_sde}
\end{equation}
with
\begin{equation}
K = \int_{\T^d} \left( I + \nabla_y \phi(y) \right)  \left( I + \nabla_y \phi(y)
\right)^T \, \mu(dy) \label{e:coeffs}
\end{equation}
and
\begin{equation}
\mu(dy) = \rho(y) dy = \frac{1}{Z} e^{-p(y)/\sigma} \, dy, \quad Z = \int_{\T^d}
e^{-p(y)/\sigma} \, dy. \label{e:gibbs_torus}
\end{equation}
The field $\phi(y)$ is the solution of the Poisson equation
\begin{equation}
- \LL_0 \phi(y) = -\nabla_y p(y), \quad \LL_0 := - \nabla_y p(y) \cdot \nabla_y + \sigma
\Delta_y, \label{e:cell}
\end{equation}
with periodic boundary conditions. The function $\rho(y)$ spans the null-space of ${\cal
L}_0^*$, the $L^2$--adjoint of $\LL_0$. The effective diffusion tensor is positive
definite and the diffusivity is always depleted \cite{Oll94}. Physically this occurs
because the homogenized process must represent the cost of traversing the many small
energy barriers present in the original multiscale problem but which are not explicitly
captured in the homogenized potential.  In Figure \ref{fig:potential} we
plot the potential $V^\eps(x,x/\eps)$, as well as the average potential $V(x)$,
illustrating this phenomenon.  In fact, the effective diffusivity $\Sigma = \sigma K$
decays exponentially fast in $\sigma$ as $\sigma \rightarrow 0$.
See \cite{CampPiatn2002} and the references therein. Thus the original
and homogenized diffusivities are exponentially different at small
temperatures.

To illustrate these facts explicitly, consider the problem in one dimension, $d = 1$. In
this case the limiting equation takes the form
\begin{equation}
 d x(t) = - A  V'(x(t)) dt + \sqrt{2 \Sigma } d \beta (t).
\label{e:lim_sde_1d}
\end{equation}
The effective coefficients are
\begin{equation}
A = \frac{\alpha L^2}{Z \widehat{Z}} \quad \mbox{and} \quad \Sigma = \frac{\sigma L^2}{Z
\widehat{Z}}, \label{e:coeffs_1d}
\end{equation}
where
\begin{equation}
\widehat{Z} = \int_{0}^L e^{p(y)/\sigma} \, dy, \quad Z = \int_{0}^L e^{-p(y)/\sigma} \,
dy. \label{e:z_1d}
\end{equation}
Notice that $L^2 \leq Z \widehat{Z}$ by the Cauchy--Schwarz inequality. This explicitly shows
that the homogenized equation in one dimension comprises motion in the average potential
$V(x)$, at a new slower time--scale contracted by $A/\alpha.$

The main results of the paper can be summarized as follows.
Assume that we are given a path $\{x^\eps(t)\}_{t \in [0,T]}$
of equation \eqref{e:xeps_V} and that we want to fit an SDE of the form
\eqref{e:lim_sde_1d} to the given data, estimating the parameters $A,\Sigma$
as $\widehat{A}, \widehat{\Sigma}$. Then the following is a loose
statement of our main results; these will be formulated precisely, and
proved, below.

\begin{theorem}
If we do not subsample, then the estimators $\widehat{A}$ and $\widehat{\Sigma}$ are
asymptotically biased -- they converge to $\alpha, \, \sigma$.
\end{theorem}
\begin{theorem}
If the sampling rate is between the two characteristic time scales of the SDE \eqref{e:main}
then the estimators $\widehat{A}$ and $\widehat{\Sigma}$ are
asymptotically unbiased -- they converge to $A, \, \Sigma$.
\end{theorem}
The rest of the paper is organized as follows. In section \ref{sec:estim} we present the
estimators that we will use. In section \ref{sec:numerics} we present various numerical
experiments illustrating the behaviour of these estimators.
In section \ref{sec:results} we state the main results of this paper,
explaining the numerical experiments from the previous section.
Section 5 contains some preliminary results that will be useful in the sequel.
Section 6 contains proof of two central propositions concerning the behaviour
of the multiscale diffusion
when observed on time--scales long compared with the fast time--scales of process, but
small compared with the slow time--scales of the process.  Section 7 is devoted to
the proofs of our theorems. Finally, section \ref{sec:conc} is devoted to some concluding
remarks.

In the sequel we use $\langle \cdot, \cdot \rangle$ to denote the standard inner--product
on $\bbR^d$ and $|\cdot|$ the induced Euclidean norm. Throughout the paper we make the
following standing assumptions on the drift vector fields:
\begin{assumptions}
\label{a:1}
The potentials $p$ and $V$ satisfy:
\begin{itemize}
\item $p(y) \in C^{\infty}_{per}(\bbT^d,\bbR^d)$;
\item $V(x) \in C^{\infty}(\R^d,\R);$
\item $|\nabla V(x_1)-\nabla V(x_2)| \le L |x_1-x_2| \quad \forall x_1,x_2 \in \R^d;$
\item $\exists a,b>0: \la -\nabla V(x),x \ra \le a-b|x|^2 \quad \forall x \in \R^d;$
\item $e^{-\frac{\alpha}{\sigma} V(x)} \in L^1(\R^d,\R^+)$.
\end{itemize}
\end{assumptions}
The third assumption will be used primarily to deduce that, by choice of
origin for $V$,
\begin{equation}
\label{e:linbnd}
|\nabla V(x)| \le L|x|.
\end{equation}
This assumption could be relaxed and replaced by a polynomial growth bound;
however this complicates the analysis without adding new insight.
Similarly it is not necessary, of course, that $V$ and $p$ are $C^{\infty}$.
The fourth condition, however, is essential:
it drives the ergodicity of the process which we use in a fundamental way in the analysis
of the drift parameter estimators; it would not, however, be fundamental for estimation
of diffusion coefficients alone. The fourth condition implies the fifth, which is simply
the requirement that the invariant measure is indeed a probability measure; we state the
two conditions separately for clarity of exposition.
%
%
%%%%%%%%%%%%%%%%%%%%%%%%%%%%%%%%%%%%%%%%%%%%%%%%%%%%%%%%%%%%%%%%%%%%%%%%%%%%%%%%%%%%%%%%%%%%%%
%
%%%%%%%%%%%%%%%%%%%%%%%%%%%%%%%%%%%%%%%%%%%%%%%%%%%%%%%%%%%%%%%%%%%%%%%%%%%%%%%%%%%%%%%%%%%%%
%
\section{The Estimators}\label{sec:estim}
In this section we describe various estimators for the parameters
arising in equation \eqref{e:lim_sde}. We assume that we are given
a path $x=\{x(t)\}_{t \in [0,T]}$, or samples from such a path,
$x=\{x_n\}_{n=0}^N$, with $x_n=x(n\delta).$ For simplicity
we aim to fit the equation in the form
\begin{equation}
 d x (t) = -A \nabla V(x (t)) dt + \sqrt{2 \Sigma } d \beta (t),
\label{e:lim_sde_sim}
\end{equation}
where $A$ and $\Sigma$ are scalars. In one dimension this reduces to the form
\eqref{e:lim_sde_1d}. Note that in general this is only the correct form for the
homogenized equation in one dimension since, typically, the average potential has a
matrix as a pre--factor, as in \eqref{e:lim_sde}. However it suffices to exemplify the
main ideas in this work, and simplifies the presentation.

The standard way to estimate the diffusion coefficient
is via the quadratic variation of the path:
\begin{equation}
\widehat{\Sigma}_{N,\delta}(x) = \frac{1}{2 N \delta d} \sum_{n = 0}^{N-1}
|x_{n+1} -  x_n|^2.
\label{e:sigma_estim_1d}
\end{equation}
A key issue in this paper is to understand how to choose $\delta$
as a function of $\epsilon$ to ensure that data generated by \eqref{e:main}
can be effectively fit to obtain the correct homogenized diffusivity
in equations such as  \eqref{e:lim_sde_sim}.

The standard way to estimate drift coefficients is via the path-space likelihood of
\eqref{e:lim_sde_sim} with respect to a pure diffusion with no drift,
namely (see, for example, \cite{BasRao80,LipShir01a})
$$L(x) \propto  \exp\{-I(x)/2\Sigma\}$$
where
$$I(x)=\int_0^T\left\{ |A \nabla V(x (t))|^2dt+2A \la \nabla V(x (t)), d x (t) \ra
\right\}.$$
Maximizing the log-likelihood then gives
the estimate $\widehat{A}$ of $A$ given by
\begin{equation}
\widehat{A}(x) = - \frac{\int_0^T \la \nabla V(x(t)), d x(t)  \ra} {\int_0^T \big| \nabla
V(x(t)) \big|^2\, dt}. \label{e:a_est}
\end{equation}
If the data is given in discrete but finely spaced increments, as often happens
in practice, then this estimator can be approximated to yield
\begin{equation}
\widehat{A}_{N,\delta}(x) = - \frac{\sum_{n = 0}^{N-1} \la \nabla V(x_n),
\left(x_{n+1} - x_n \right)\ra}{\sum_{n=0}^{N-1} \left|\nabla V(x_n) \right|^2
\delta}.
\label{e:alpha_estim_1d}
\end{equation}
A key issue in this paper is to understand how to chose $\delta$ as a function of
$\epsilon$ to ensure that data generated by \eqref{e:main} can be effectively fit to
obtain the correct homogenized drift coefficients in equations such as
\eqref{e:lim_sde_sim}, via the estimator \eqref{e:alpha_estim_1d}.

The gradient structure of the SDE \eqref{e:lim_sde_sim} can be used to obtain a second
estimator for the drift coefficients. This second estimator, which we now derive, is of
interest for two different reasons: firstly it may be useful in practice as it may lead
to smaller variance in estimators; secondly it highlights the fact that working out how
to sample the data to obtain the correct estimation of the diffusion coefficient alone
will lead to correct estimation of the drift parameters, at least for the
class of gradient--structure SDEs that we consider in this paper.
The second estimator requires
the input of an estimator $\widehat\Sigma$ for the diffusion coefficient and is
\begin{equation}
\tilde{A}(x)=\widehat{\Sigma}\frac{\frac{1}{T}\int_0^T \Delta V(x(t)) \, dt  }
{\frac{1}{T} \int_0^T |\nabla V(x(t))|^2 \, dt}. \label{eq:alpha2}
\end{equation}
Approximating to allow for the input of discrete--time data gives
\begin{equation}
\tilde{A}_{N,\delta}(x) =  \widehat{\Sigma}\frac{\sum_{n = 0}^{N-1} \Delta V(x_n) \delta}
{\sum_{n=0}^{N-1} \left|\nabla V(x_n) \right|^2 \delta}.
\label{e:alpha_estim_1d2}
\end{equation}
The following result shows that $\tilde{A}(x)$ is a natural approximation to $\widehat
A(x).$
\begin{prop}
Let $x=\{x(t)\}_{t \in [0,T]}$ satisfy \eqref{e:lim_sde_sim}. If $\widehat{\Sigma}=\Sigma$
then the estimator $\tilde A(x)$ is asymptotically equivalent to the maximum likelihood
estimator $\widehat{A}$:
$$
\lim_{T \rightarrow \infty} \tilde{A}(x) = \widehat{A}(x),\, a.s.
$$
\end{prop}
\proof We apply the It\^{o} formula to $V(x(t))$ for $x(t)$ solving \eqref{e:lim_sde_sim}
and use formula \eqref{e:a_est} to obtain
\begin{eqnarray*}
\widehat{A}(x) & = &
\frac{V(x(0)) - V(x(T)) + \Sigma \int_0^T \Delta V(x(t))
\, dt  }{\int_0^T |\nabla V(x)|^2 \, dt}
\\ & = & \frac{(V(x(0)) - V(x(T)))}{\int_0^T |\nabla V(x)|^2 \,
dt} + \frac{\frac{1}{T}\Sigma\int_0^T \Delta V(x(t)) \, dt  } {\frac{1}{T}
\int_0^T |\nabla V(x)|^2 \, dt}
\\ & = & \frac{\frac{1}{T} (V(x(0)) - V(x(T)))}{\frac{1}{T}\int_0^T |\nabla V(x)|^2 \,
dt} + \tilde{A}(x).
\end{eqnarray*}
Under the Assumptions \ref{a:1} it follows from \cite{Mao97} that
$$
\lim_{T \rightarrow 0} \frac{\frac{1}{T} (V(x(0)) - V(x(T))}{ \int_0^T
|\nabla V(x(t))|^2 \, dt} = 0,\, a.s.
$$
The result follows.
\qed
%
%
%%%%%%%%%%%%%%%%%%%%%%%%%%%%%%%%%%%%%%%%%%%%%%%%%%%%%%%%%%%%%%%%%%%%%%%%%%%%%%%%
%
%                          NUMERICAL RESULTS
%
%%%%%%%%%%%%%%%%%%%%%%%%%%%%%%%%%%%%%%%%%%%%%%%%%%%%%%%%%%%%%%%%%%%%%%%%
%
\section{Numerical Results}
\label{sec:numerics} In all cases we solve the multiscale SDE \eqref{e:main} using the
Euler--Marayama scheme \cite{KlPl92} for a single realization of the noise, with a
time--step $\Delta t$ sufficiently small so that the error due to the discretization is
negligible; this requires that the time--step is small compared with $\eps^2,$ the
fastest scale in the problem. We also employ a sufficiently long time interval so that
the invariant measure is well sampled by the single path. Since the convergence to the invariant measure
is uniform in $\eps \to 0$, this is not prohibitive. We then use the data generated from
the multiscale process as input to the estimators for the homogenized diffusion
\eqref{e:lim_sde}. We present numerical results for three model problems: a one
dimensional monomial potential of even degree, a one dimensional bistable potential and a
two dimensional quadratic potential. In all three cases we perturb the large--scale part
of the potential $V$ by small--scale fast oscillations, usually in the form of a cosine
potential $p$.

We present two types of numerical results. Note that $\delta$, the time interval between
two consecutive observations, is the inverse sampling rate. In the first we use
$\delta=\Delta t$ as the time interval between two consecutive observations in the
estimators. In the second we subsample the data, using $\delta > \Delta t$ and study how
the estimated coefficients behave as a function of the subsampling. We use the data
generated from our simulation in the estimators \eqref{e:alpha_estim_1d} and
\eqref{e:alpha_estim_1d2} to estimate the drift coefficient and in \eqref{e:sigma_estim_1d}
to estimate the diffusion coefficient of \eqref{e:lim_sde_1d}. For the most part we work
in one dimension and fit a single drift and diffusion parameter so that \eqref{e:lim_sde}
becomes \eqref{e:lim_sde_1d}. When we work in more than one dimension, or
estimate more than just a single drift or diffusion parameter,
we use natural generalizations of the estimators defined in the previous
section.

Let us summarize the main conclusions that can be drawn from the numerical experiments;
recall that $\Delta t \ll \eps^2.$ First, if we choose $\delta = \Delta t$, that is, if
we don't subsample, then the resulting estimators do not generate the correct estimates
of the homogenized coefficients. If, on the other hand, we subsample with $\eps^2 \ll
\delta \ll \mathcal{O} (1),$ then  the estimators generate the values of the parameters
of the homogenized equation. Furthermore, there is an optimal sampling rate: there exists
a $\delta^*$ which minimizes the distance between the homogenized value of the parameter
and the value generated by the estimator. The optimal sampling rate depends sensitively
on $\sigma$. It is also of interest that, in higher dimensions, the optimal sampling rate
can be different for different parameters.

The above observations appear to hold independently of the detailed
form of the large--scale part of the potential $V$ (provided, of course,
that it satisfies appropriate convexity conditions).
In addition, the performance of the estimators seems to be the same
irrespective of the dimension of the problem.

Another interesting observation is that the second estimator for the drift coefficient
\eqref{e:alpha_estim_1d2} performs at least as well as the maximum
likelihood estimator \eqref{e:alpha_estim_1d}, and in some instances
outperformas it.
\subsection{Failure Without Subsampling}
\begin{figure}
\centerline{
\begin{tabular}{c@{\hspace{2pc}}c}
\includegraphics[width=2.7in, height = 2.7in]{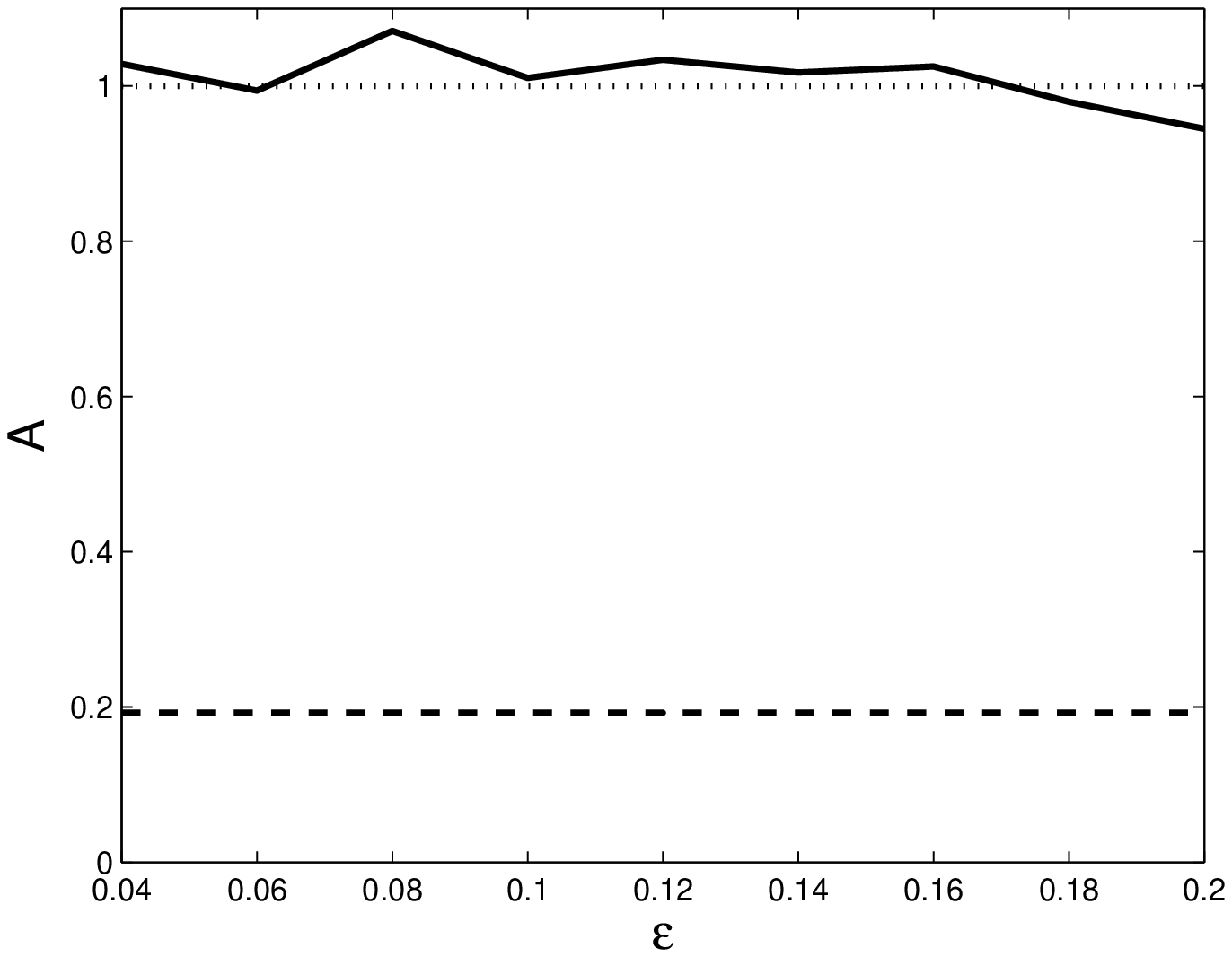} &
\includegraphics[width=2.7in, height = 2.7in]{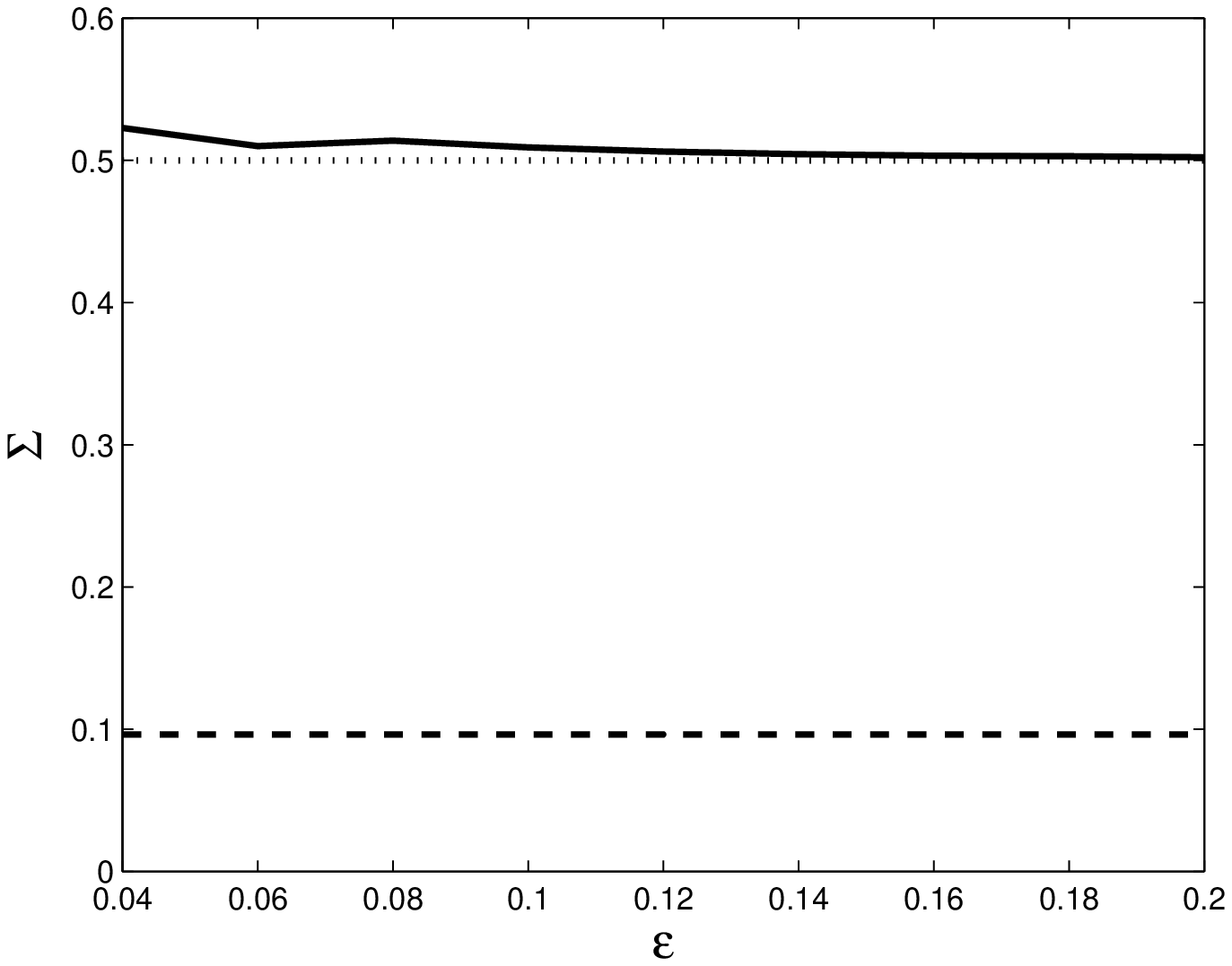} \\
 a.~~  $\widehat{A}$  & b.~~ $\widehat{\Sigma}$
\end{tabular}}
\begin{center}
\caption{ Estimation of the drift and diffusion coefficients vs $\eps$ for the potential \eqref{e:ou}.
Solid line: estimated coefficient. Dashed line: homogenized coefficient. Dotted line:
unhomogenized coefficient.}
\label{fig:vs_eps_no_subsam}
\end{center}
\end{figure}
\begin{figure}
\centerline{
\begin{tabular}{c@{\hspace{2pc}}c}
\includegraphics[width=2.7in, height = 2.7in]{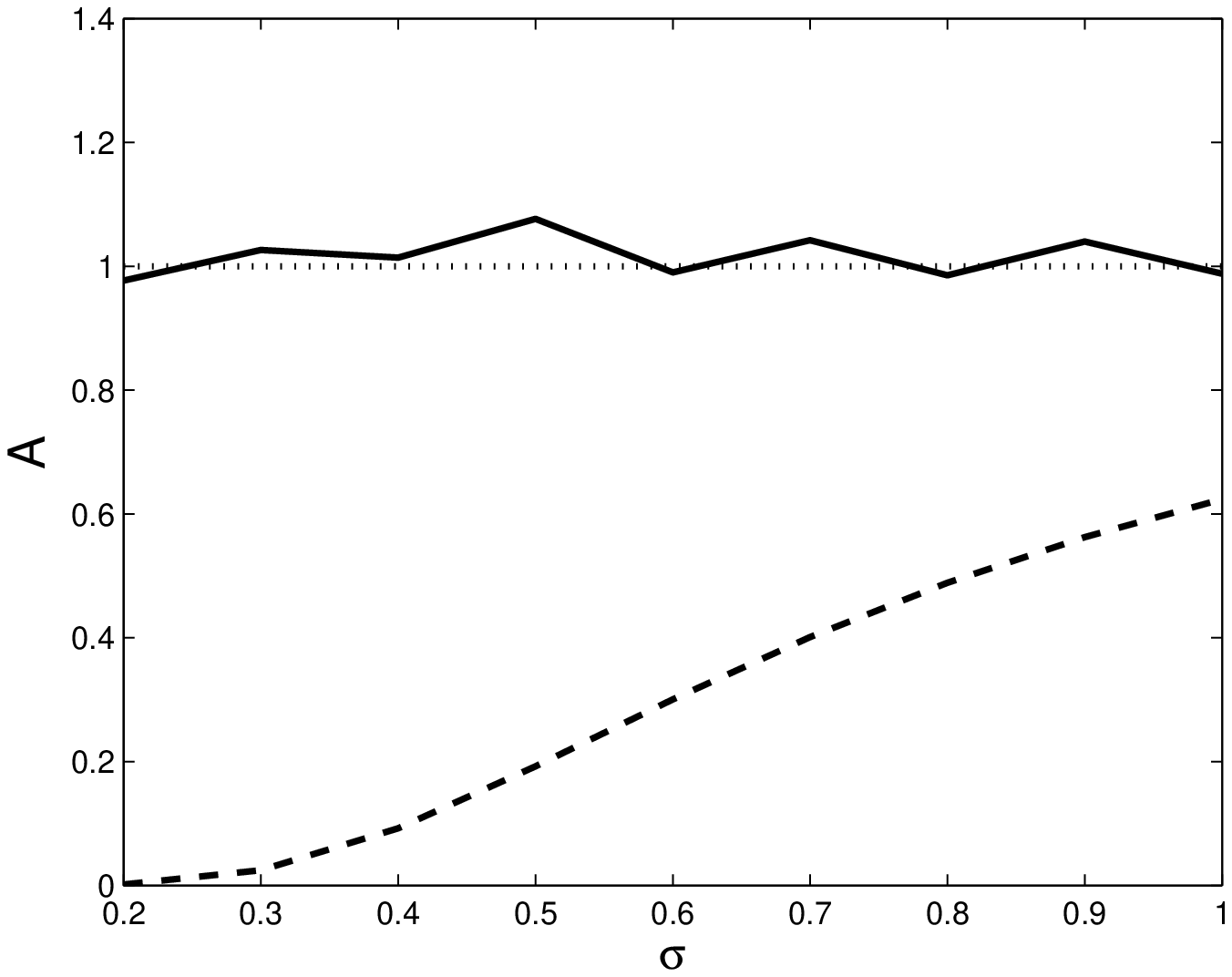}
 & \includegraphics[width=2.7in, height = 2.7in]
 {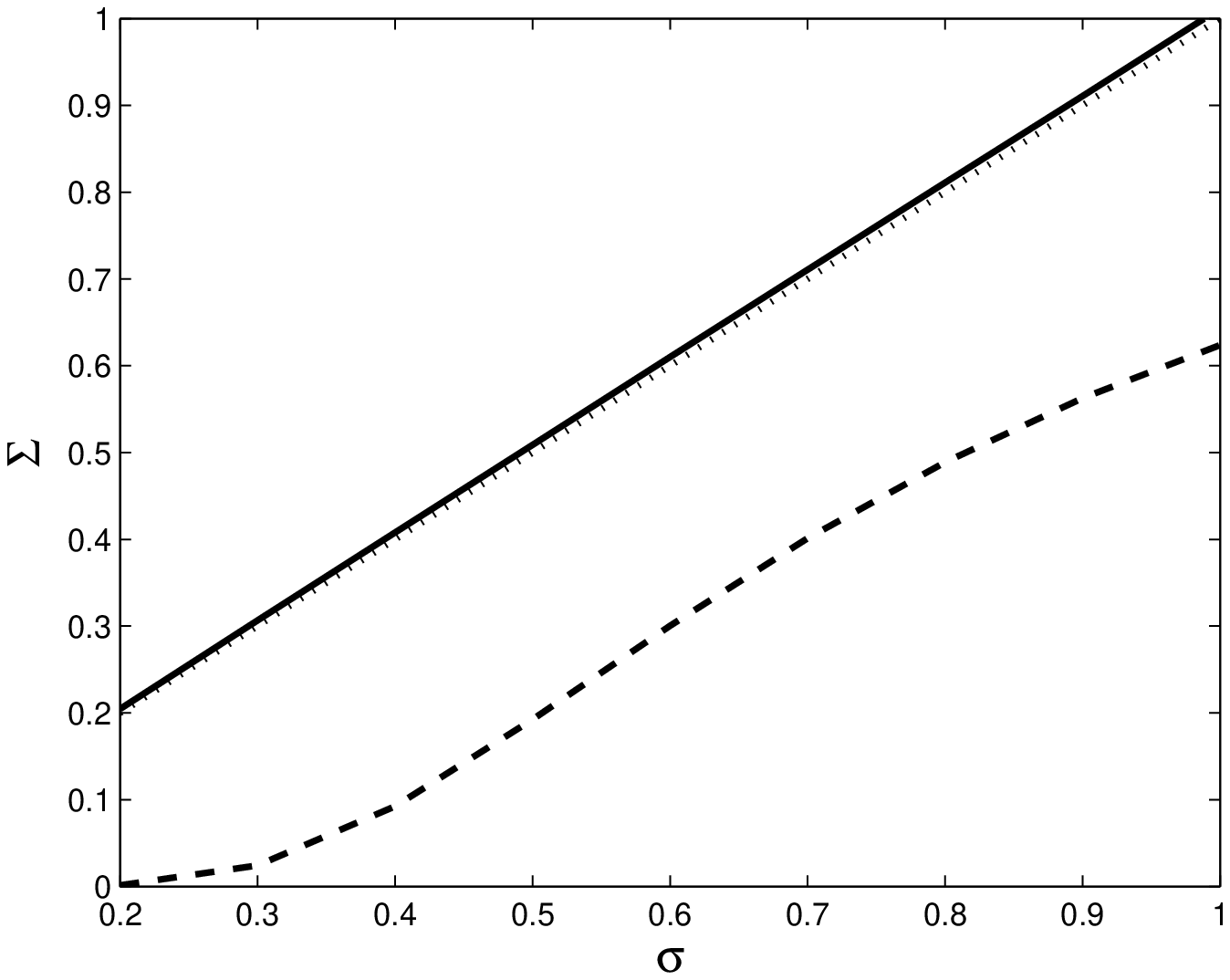} \\
 a.~~  $\widehat{A}$  & b.~~ $\widehat{\Sigma}$
\end{tabular}}
\begin{center}
\caption{Estimation of the drift and diffusion coefficients vs $\sigma$ for the potential
\eqref{e:ou} with $\eps = 0.1$. Solid line: estimated coefficient. Dashed line:
homogenized coefficient. Dotted line: unhomogenized coefficient.}
\label{fig:vs_sigma_no_subsam}
\end{center}
\end{figure}
In this section we study the estimators
$\widehat{A}$ and $\widehat{\Sigma}$
when the data is given from the solution of equation \eqref{e:xeps_V}
with $\eps \ll 1$ and $\Delta t=\delta$ -- no subsampling is used.
We use the potential
\begin{equation}\label{e:ou}
V(x) = \frac{1}{2}\alpha x^2
\end{equation}
The small--scale part of the potential is
\begin{equation}\label{e:cos}
p(y) =  \cos ( y ).
\end{equation}
In Figure \ref{fig:vs_eps_no_subsam} we plot the estimators $\widehat{A}$ and
$\widehat{\Sigma}$ for various values of $\eps$. For comparison we also plot the
homogenized coefficients $A$ and $\Sigma$ and the unhomogenized coefficients $\alpha$ and
$\sigma$. We observe that the estimators always give us the coefficients $\alpha$ and
$\sigma$ of the original SDE \eqref{e:xeps_V}. In particular, the performance of the
estimators does not improve as $\eps \rightarrow 0$. In Figure
\ref{fig:vs_sigma_no_subsam} we plot the estimators for various values of the diffusion
coefficient $\sigma$. We notice that the estimators give the values of the coefficients
$\alpha$ and $\sigma$, for all values of $\sigma$. Since the homogenized coefficients
decay to $0$ exponentially fast in $\sigma$, the results of Figure
\ref{fig:vs_sigma_no_subsam} indicate that the estimators give exponentially wrong
results when $\sigma \ll 1$.

These results indicate the need to subsample -- i.e. to choose $\delta$
appropriately as a function of $\epsilon$.
%
%
%%%%%%%%%%%%%%%%%%%%%%%%%%%%%%%%%%%%%%%%%%%%%%%%%%%%%%%%%%%%%%%%%%%%%%%%%%%%%%%%%%%%%%%%%%%
\subsection{Success With Subsampling}
Now, rather than using all the data that were generated from the solution of equation
\eqref{e:main} we use only a fraction of them. We choose $\delta$ in the estimators
\eqref{e:sigma_estim_1d}, \eqref{e:alpha_estim_1d} and \eqref{e:alpha_estim_1d2} as
follows:
$$
\Delta t_{sam}=\delta = 2^k \Delta t, \quad k=0, \, 1, \, 2, \dots,
$$
and we study the performance of the estimators as a function of the sampling rate. We
investigate this issue for three different model problems.
\subsubsection{OU Processes in 1D}
\begin{figure}
\centerline{
\begin{tabular}{c@{\hspace{2pc}}c}
\includegraphics[width=2.7in, height = 2.7in]{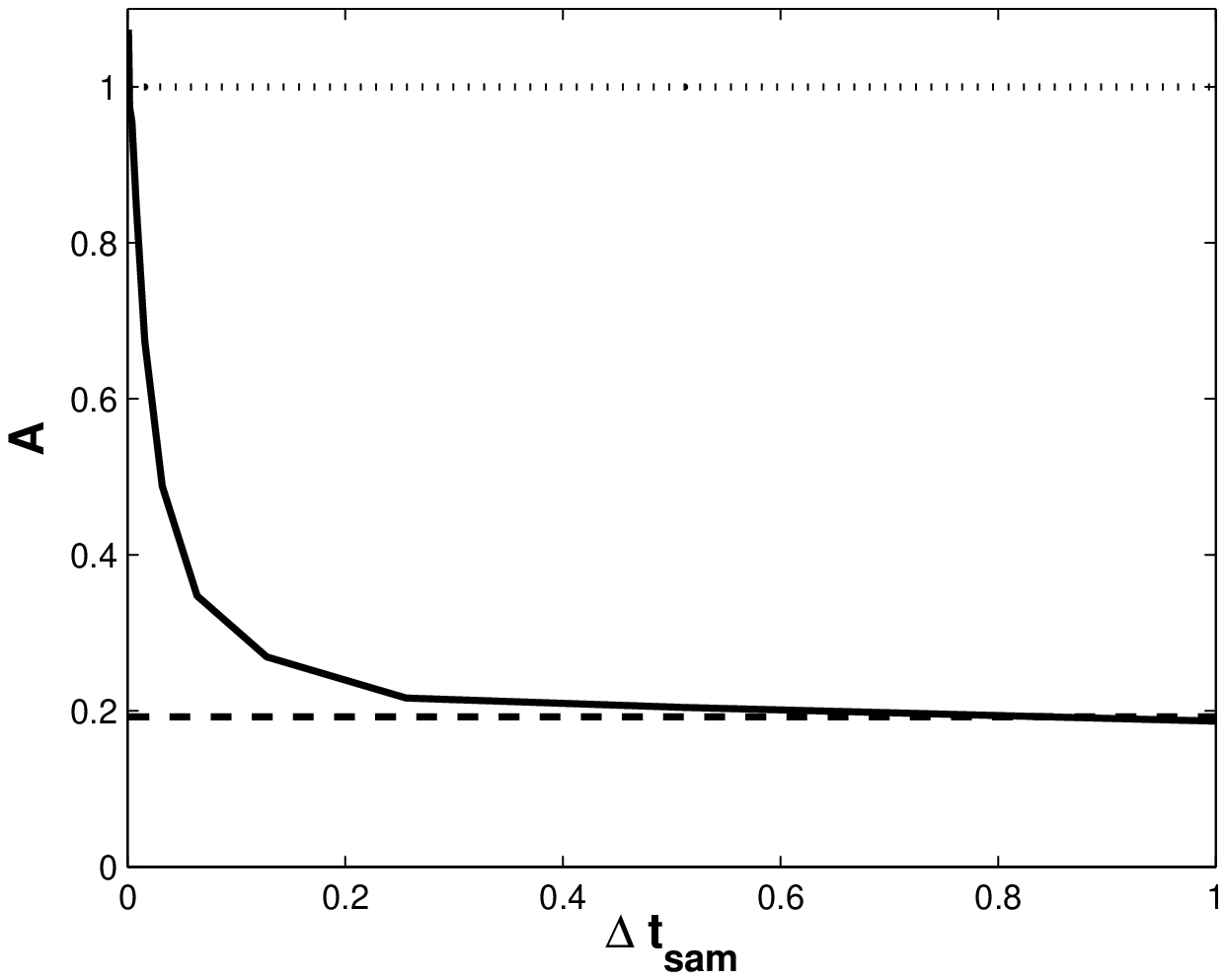}
& \includegraphics[width=2.7in, height = 2.7in]
{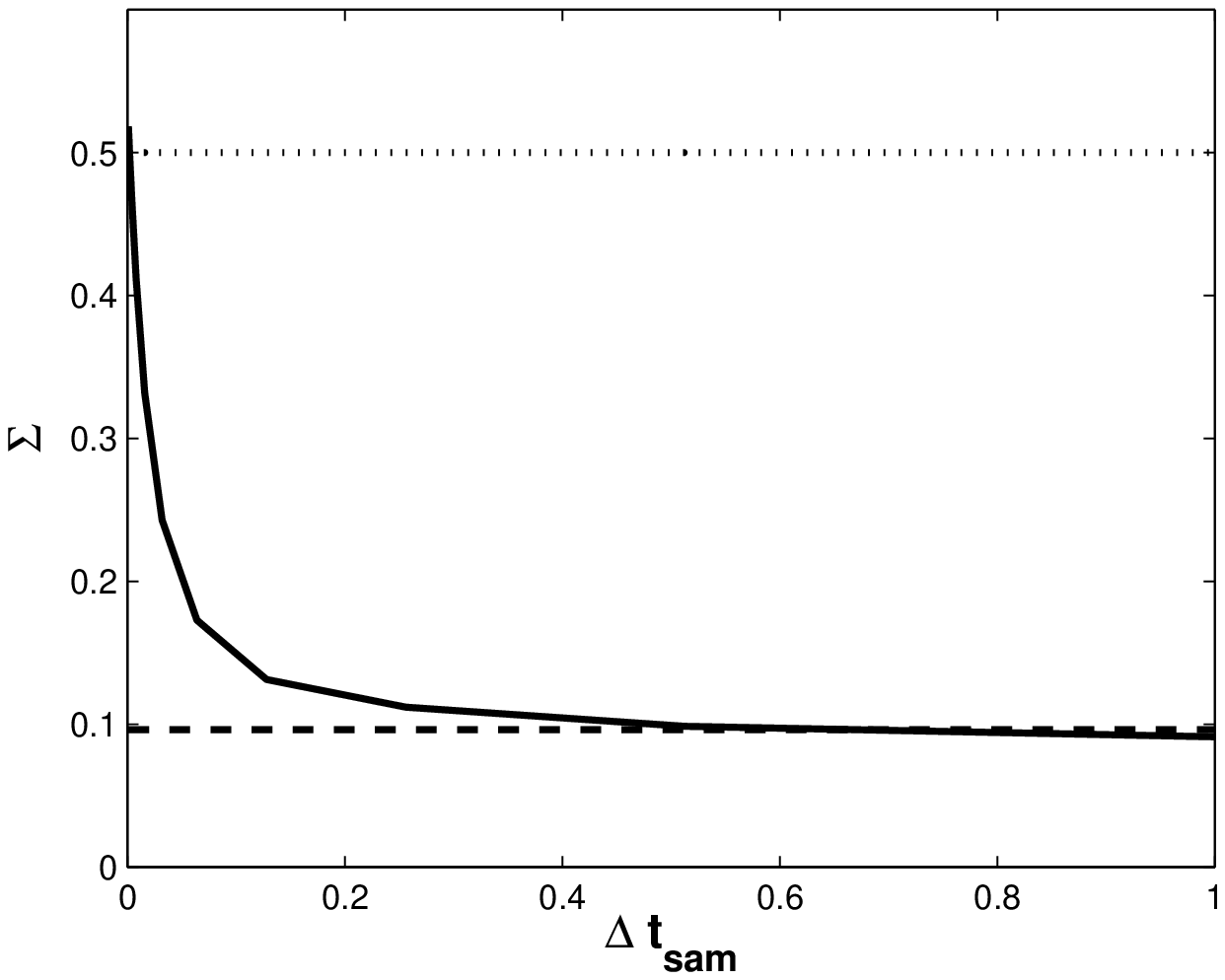} \\
a.~~  $\widehat{A}$  & b.~~ $ \widehat{\Sigma}$
\end{tabular}}
\begin{center}
\caption{Estimation of the drift and diffusion coefficients vs $\Delta t_{sam}$
for the potential \eqref{e:ou} with $\eps = 0.1$.  Solid line: estimated coefficient.
Dashed line: homogenized coefficient. Dotted line: unhomogenized coefficient.}
\label{fig:ou_sig05}
\end{center}
\end{figure}
\begin{figure}[t]
\centerline{
\begin{tabular}{c@{\hspace{2pc}}c}
\includegraphics[width=2.7in, height = 2.7in]{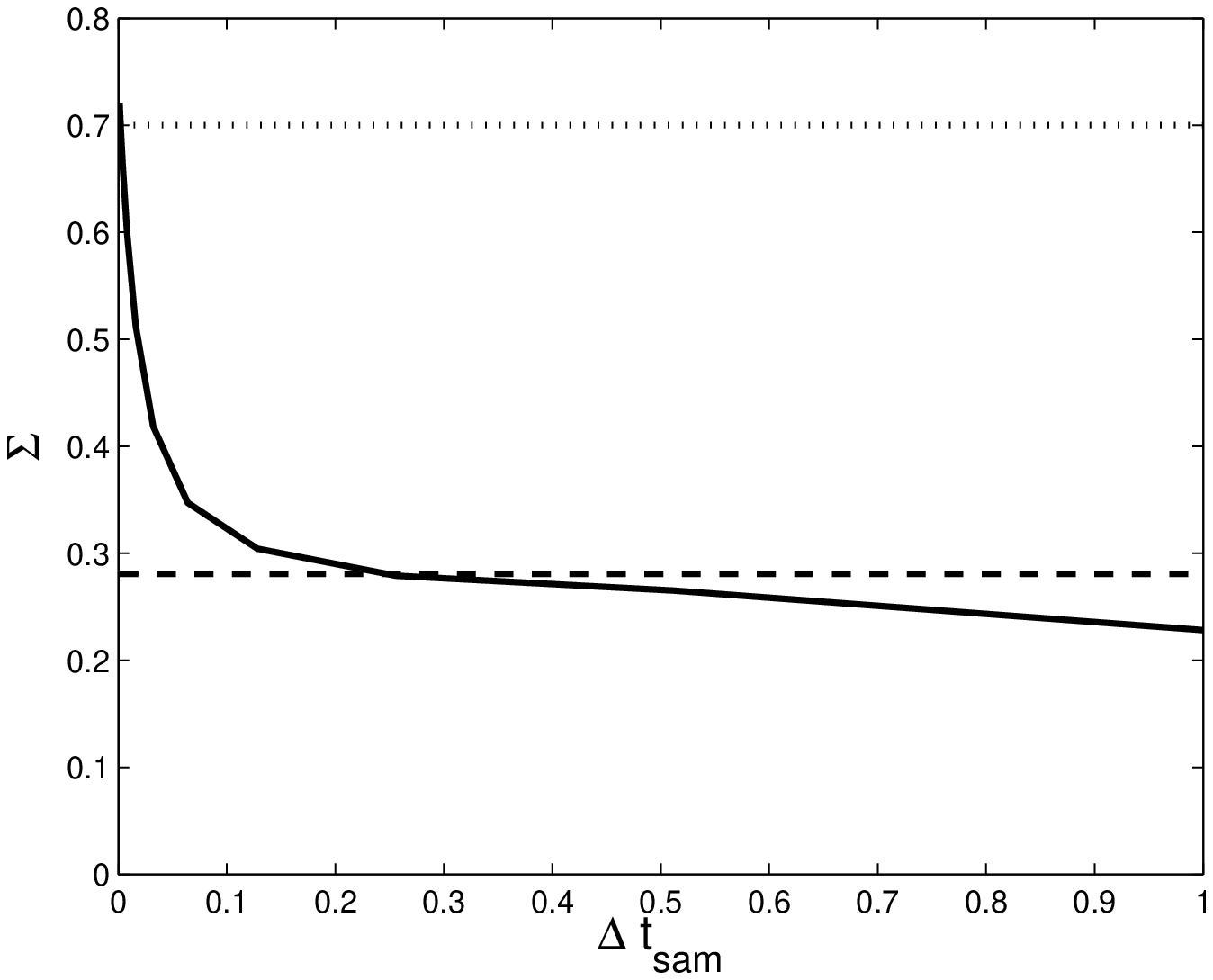}
& \includegraphics[width=2.7in, height = 2.7in]{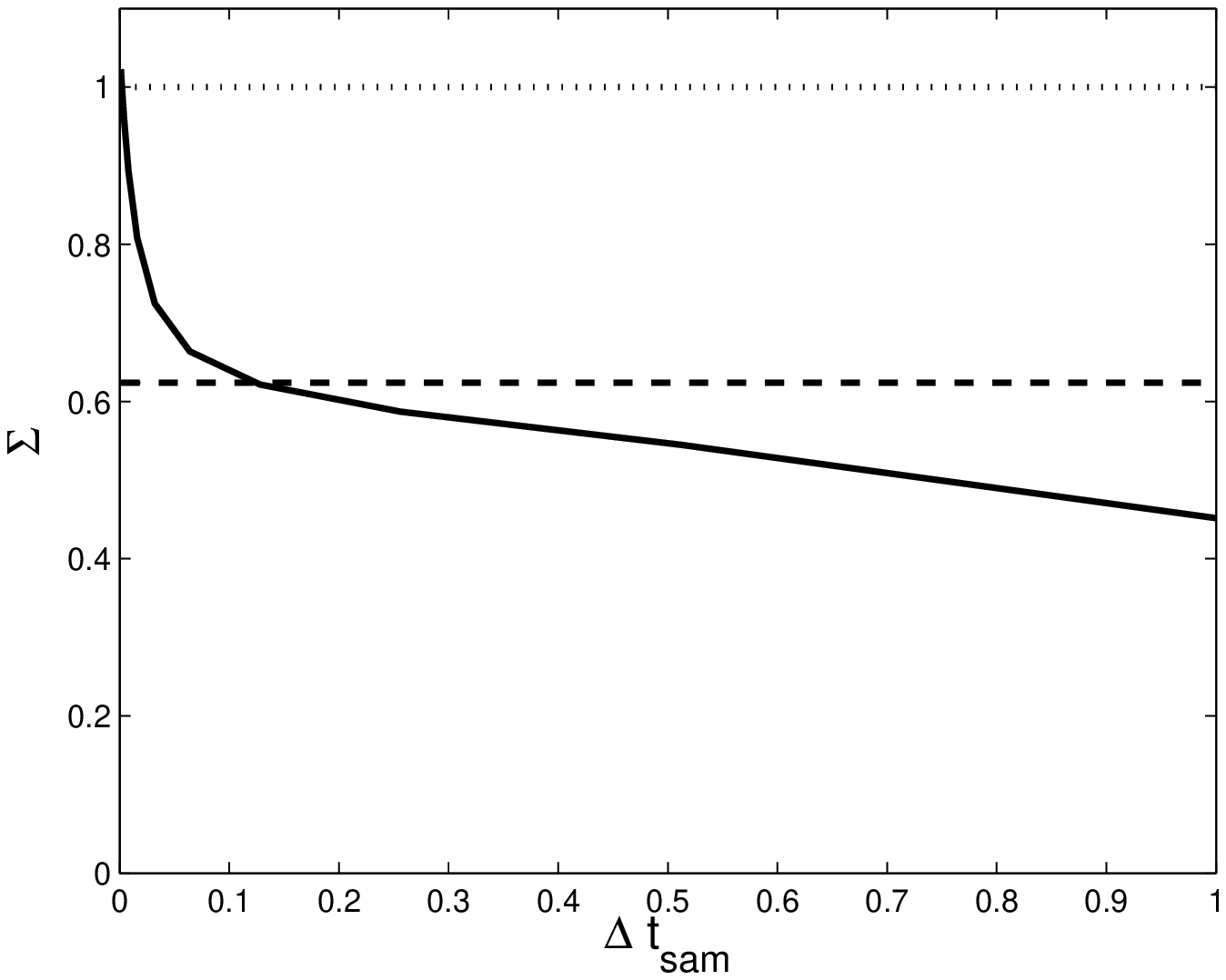}
\\  a.~~  $\sigma=0.7$  & b.~~ $\sigma = 1.0$
\end{tabular}}
\begin{center}
\caption{Estimation of the diffusion coefficient vs $\Delta t_{sam}$ for the potential
\eqref{e:ou} with $\eps = 0.1$, for two different values of $\sigma$.  Solid line: estimated
coefficient. Dashed line: homogenized coefficient. Dotted line: unhomogenized coefficient.}
\label{fig:ou_sig_07_1}
\end{center}
\end{figure}
\begin{figure}[t]
\centerline{
\begin{tabular}{c@{\hspace{2pc}}c}
\includegraphics[width=2.7in, height = 2.7in]
{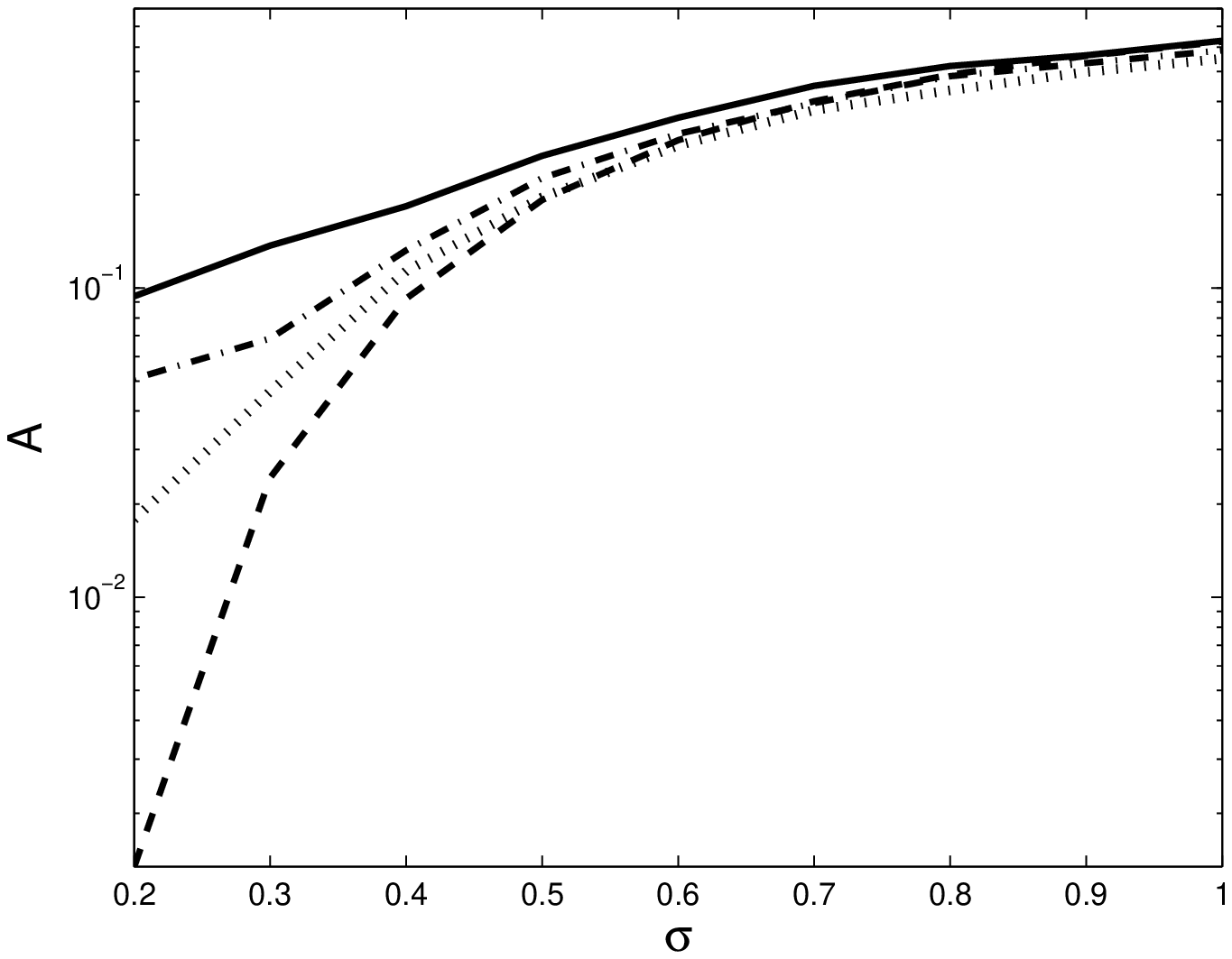} &
\includegraphics[width=2.7in, height = 2.7in]
{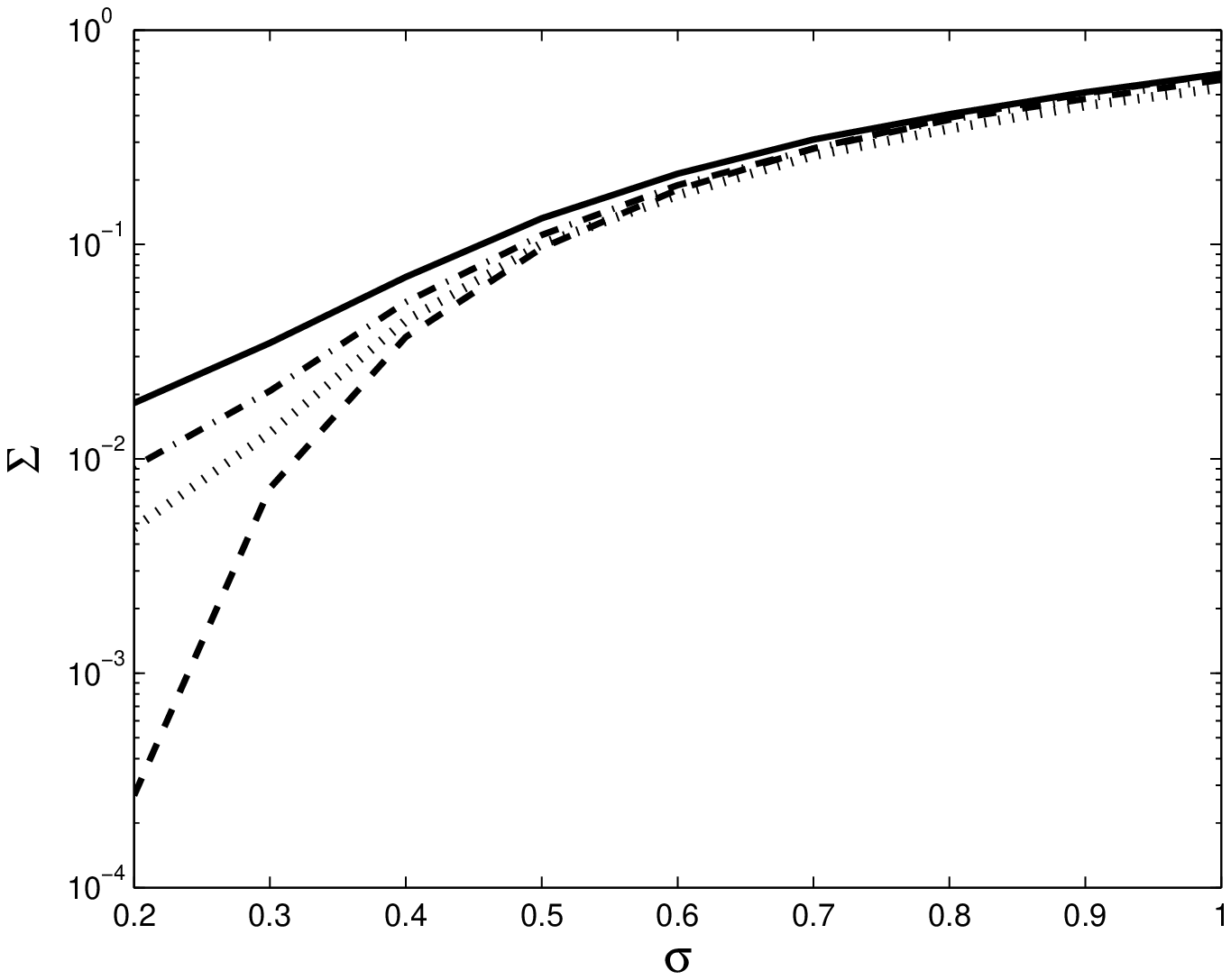} \\
 a.~~  $\widehat{A}$  & b.~~ $\widehat{\Sigma}$
\end{tabular}}
\begin{center}
\caption{Estimation of the drift and diffusion coefficient vs $\sigma$ for the potential
\eqref{e:ou} with $\eps = 0.1, \, \alpha = 1.0$, for three different sampling rates.
Solid line: $\Delta t_{sam} = 0.128$. Dash--dotted line: $\Delta t_{sam} = 0.256$. Dotted
line: $\Delta t_{sam} = 0.512$. Dashed line: homogenized coefficient. }
\label{fig:ou_vs_sig_sam}
\end{center}
\end{figure}
We study the problem in one dimension with the large--scale part of the potential given
by \eqref{e:ou} and with the fluctuating part being the cosine potential \eqref{e:cos}.
The two estimators $\widehat{A}$ and $\tilde{A}$ for the drift coefficient produce almost
identical results and we only present results for the maximum likelihood estimator
$\widehat{A}$. In Figure \ref{fig:ou_sig05} we present the estimated values of the drift
and diffusion coefficients as a function of the inverse sampling rate $\delta = \Delta
t_{sam}$ when $\eps = 0.1, \, \alpha = 1.0, \, \sigma = 0.5$. We observe that, provided
that we subsample at an appropriate rate, we are able to estimate the parameters of the
homogenized equation correctly. Notice also that the estimators for the drift and the
diffusion coefficient show very similar dependence on the sampling rate. This is in
accordance with our theoretical results; see Theorem \ref{prop:drift_estim_2}.

In Figure \ref{fig:ou_sig_07_1} we plot $\widehat{\Sigma}$ as a function of the sampling
rate for two different values of $\sigma$. We observe that the estimator of the diffusion
coefficient is a decreasing function of the sampling rate, as expected. In addition to
this, there is a well defined optimal sampling rate, which depends sensitively on
$\sigma$. In particular the optimal $\delta$ is a decreasing function of $\sigma$. This
is to be expected, since when $\sigma \gg 1$ the process $x^\eps(t)$ loses its multiscale
character and becomes effectively a standard Brownian motion. Consequently, when $\sigma$
is sufficiently large, the optimal $\delta$ becomes $\Delta t$, the integration time
step. Notice furthermore that the slope of the $\widehat{\Sigma}-\delta$
curve depends on $\sigma$.

In Figure \ref{fig:ou_vs_sig_sam} we plot the estimators of the drift and diffusion
coefficients versus $\sigma$, for three different sampling rates. For comparison we also
plot the homogenized coefficients. We observe that all three sampling rates lead to
reasonably accurate estimates for $A$ and $\Sigma$, when $\sigma$ is not too small. On
the other hand, the estimators become less accurate as $\sigma \rightarrow 0$. This is
also to be expected: when $\sigma \ll 1$, the accurate simulation of \eqref{e:main}
requires a very small time step; moreover, the equation has to be solved over a very long
time interval in order for the invariant measure of the process to be well represented.
Hence, our hypothesis that the errors due to discretization and finite time of
integration are small, is not valid. In addition, as $\sigma$ tends to $0$, the optimal sampling
rate increases, and becomes much larger than the coarser sampling rate that we use in the
simulations.

\begin{figure}[t]
\centerline{
\begin{tabular}{c@{\hspace{2pc}}c}
\includegraphics[width=2.7in, height = 2.7in]
{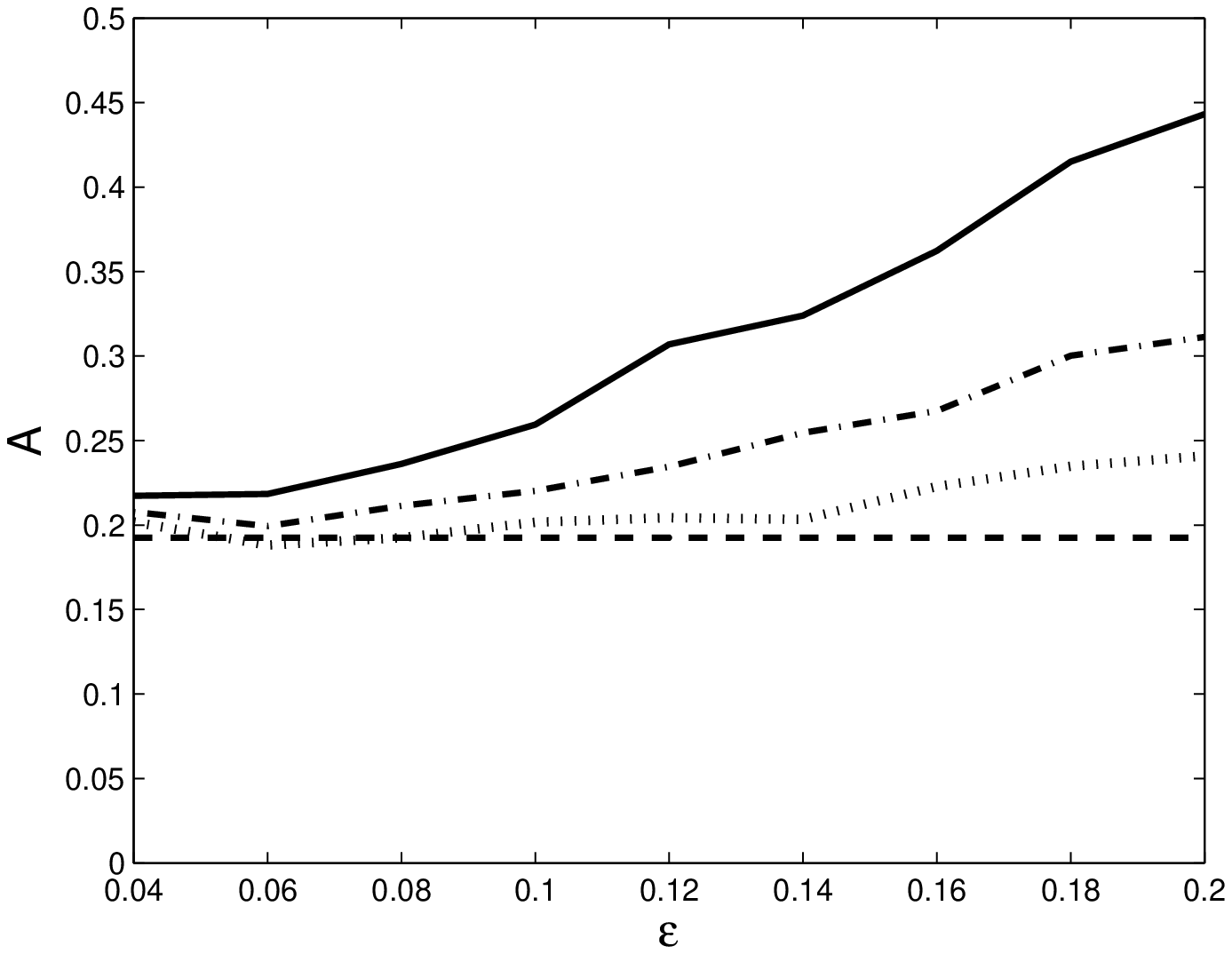} &
\includegraphics[width=2.7in, height = 2.7in]
{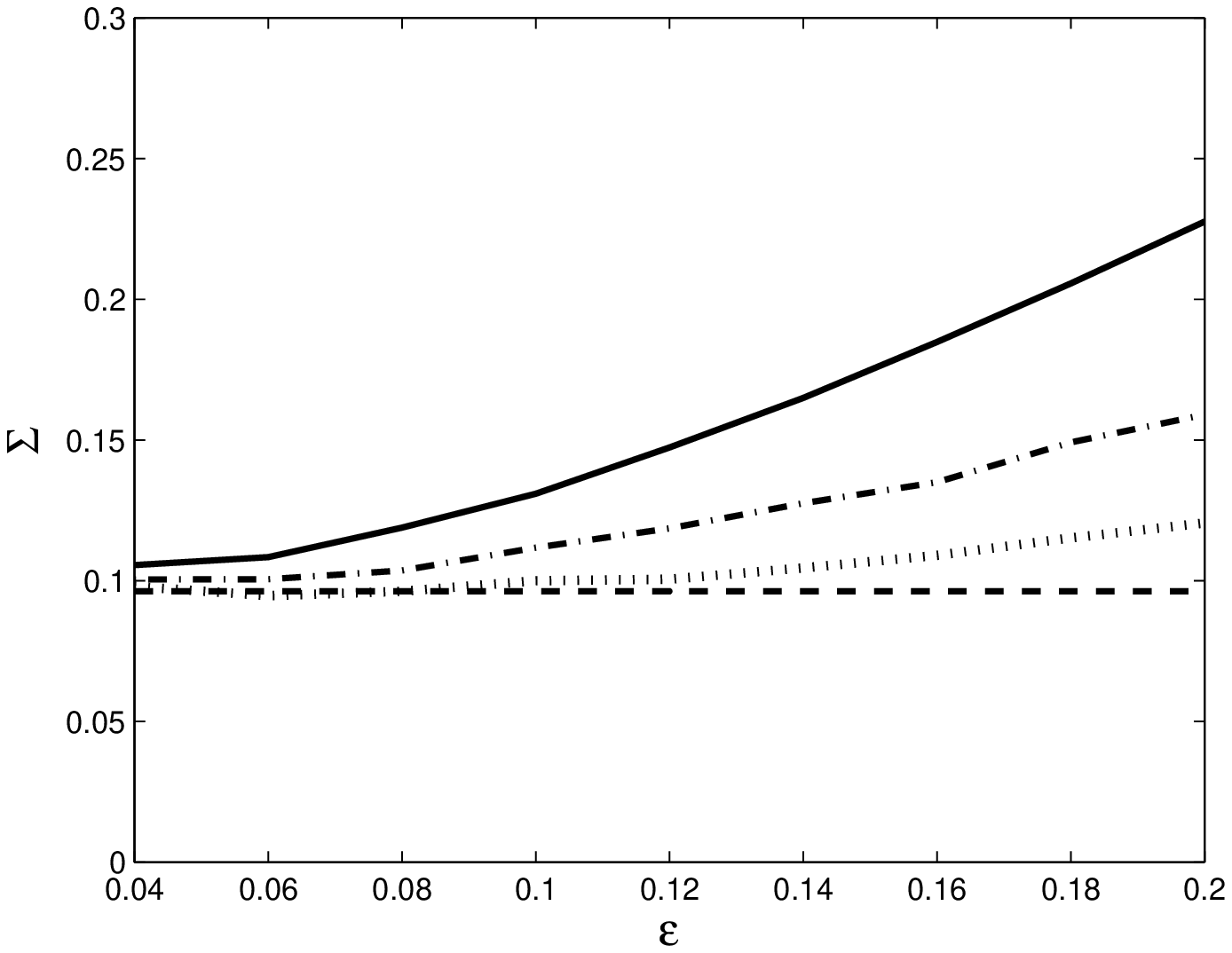} \\
 a.~~  $\widehat{A}$  & b.~~ $\widehat{\Sigma}$
\end{tabular}}
\begin{center}
\caption{Estimation of the drift and diffusion coefficient vs $\eps$ for the potential
\eqref{e:ou} with $\alpha = 1.0, \, \sigma = 0.5$, for three different sampling rates.
Solid line: $\Delta t_{sam} = 0.128$. Dash--dotted line: $\Delta t_{sam} = 0.256$. Dotted
line: $\Delta t_{sam} = 0.512$. Dashed line: homogenized coefficient.}
\label{fig:ou_vs_eps_sam}
\end{center}
\end{figure}

In Figure \ref{fig:ou_vs_eps_sam} we plot the estimators versus $\eps$, for three
different values of the sampling rate. As expected, the deviation of the estimated values
of the drift and diffusion coefficients from the homogenized values is an increasing
function of $\epsilon$. On the other hand, the optimal sampling rate does not appear to
depend sensitively on $\eps$: it is always the same sampling rate that minimizes the
distance between the estimated coefficient and the homogenized one, for all values of
$\eps$.
%
%%%%%%%%%%%%%%%%%%%%%%%%%%%%%%%%%%%%%%%%%%%%%%%%%%%%%%%%%%%%%%%%%%%%%%%%%%%%%%%%%%%%%%%%%%
%
\subsubsection{A Bistable Potential}
We consider equation \eqref{e:main} in one dimension with a mean potential of
the bistable form
\begin{equation}\label{e:pot_bistable}
V(x; \alpha, \beta) = - \frac{1}{2} \alpha x^2 + \frac{1}{4} \beta x^4.
\end{equation}
The fluctuating part of the potential is given by  \eqref{e:cos}. The homogenized
equation is
\begin{equation}\label{e:hom_bistable}
d X(t) = ( A X(t) - B X(t)^3 ) dt + \sqrt{2 \Sigma} d \beta(t),
\end{equation}
where the homogenized coefficients are given by
$$
A = \alpha K, \quad B = \beta K, \quad \Sigma = \sigma K,
\quad K = \frac{4 \pi^2}{Z \widehat{Z}},
$$
where $Z$ and $\widehat{Z}$ are given by \eqref{e:z_1d} with $L = 2 \pi$ and $p(y) = \cos(y)$.
We will estimate the diffusion coefficient using formula \eqref{e:sigma_estim_1d} with $d
= 1$. For the two parameters of the drift we use generalizations of the maximum
likelihood estimator $\widehat{A}$.

\begin{figure}[t]
\centerline{
\begin{tabular}{c@{\hspace{2pc}}c}
\includegraphics[width=2.7in, height = 2.7in]{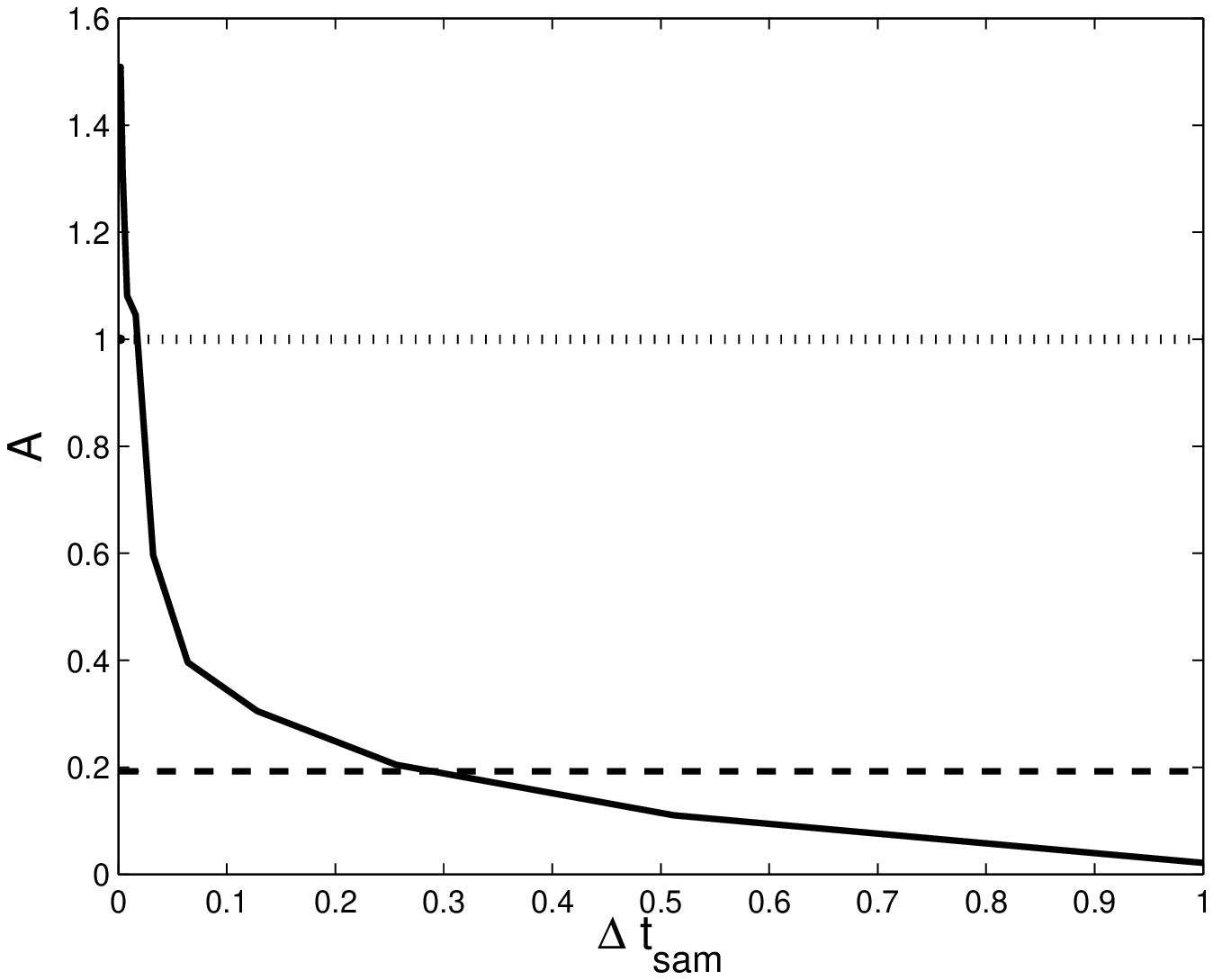} &
\includegraphics[width=2.7in, height = 2.7in]{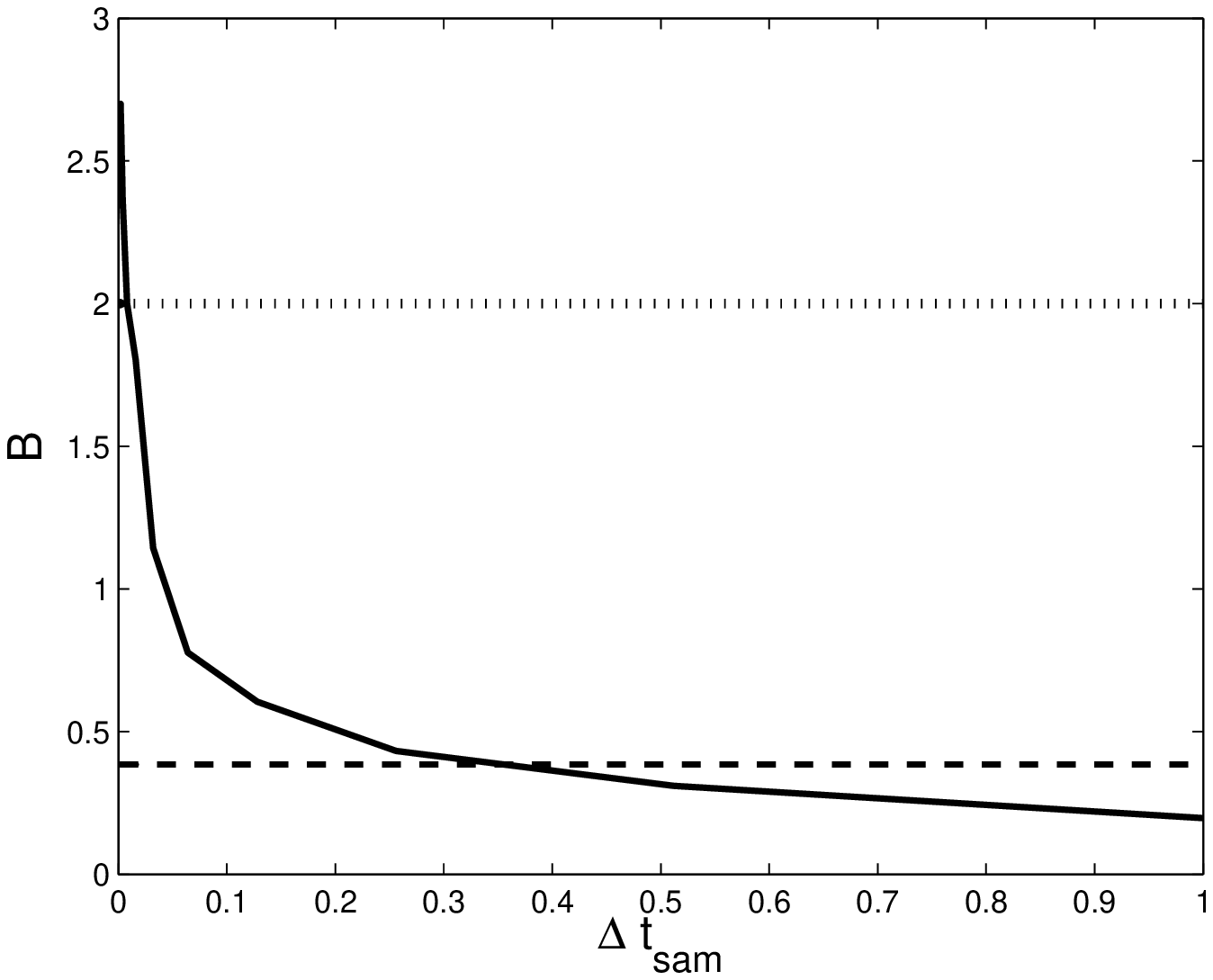} \\
a.~~  $\widehat{A}$ vs $\Delta t_{sam}$  & b.~~ $ \widehat{B}$ vs $\Delta t_{sam}$
\end{tabular}}
\begin{center}
\caption{ Estimation of the parameters of the bistable potential \eqref{e:pot_bistable}
as a function of the sampling rate for $\sigma = 0.5, \,\eps = 0.1$. Solid line:
estimated coefficient. Dashed line: homogenized coefficient. Dotted line: unhomogenized
coefficient.} \label{fig:bistable:A_B_05}
\end{center}
\end{figure}

\begin{figure}[t]
\centerline{
\begin{tabular}{c@{\hspace{2pc}}c}
\includegraphics[width=2.7in, height = 2.7in]{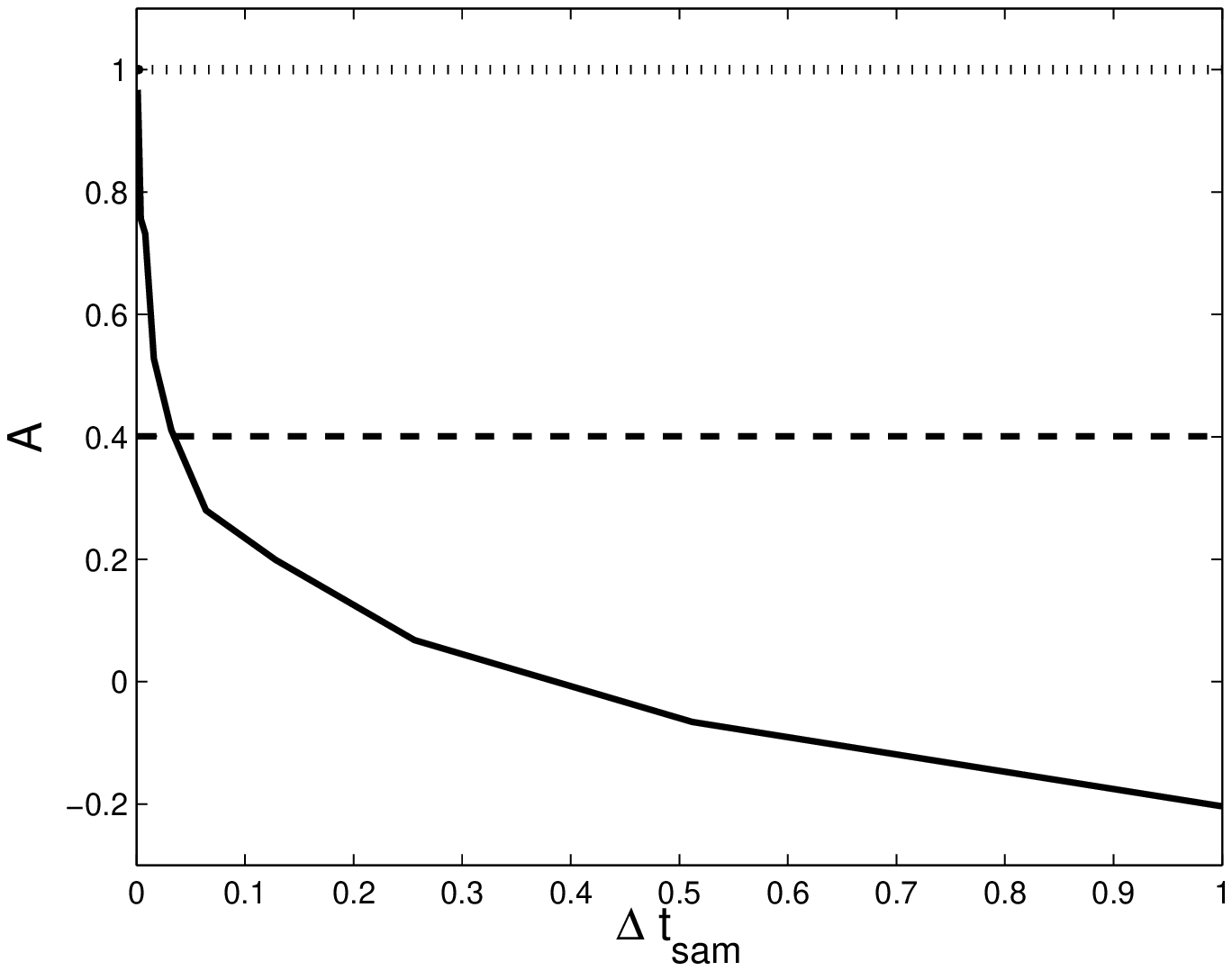} &
\includegraphics[width=2.7in, height = 2.7in]{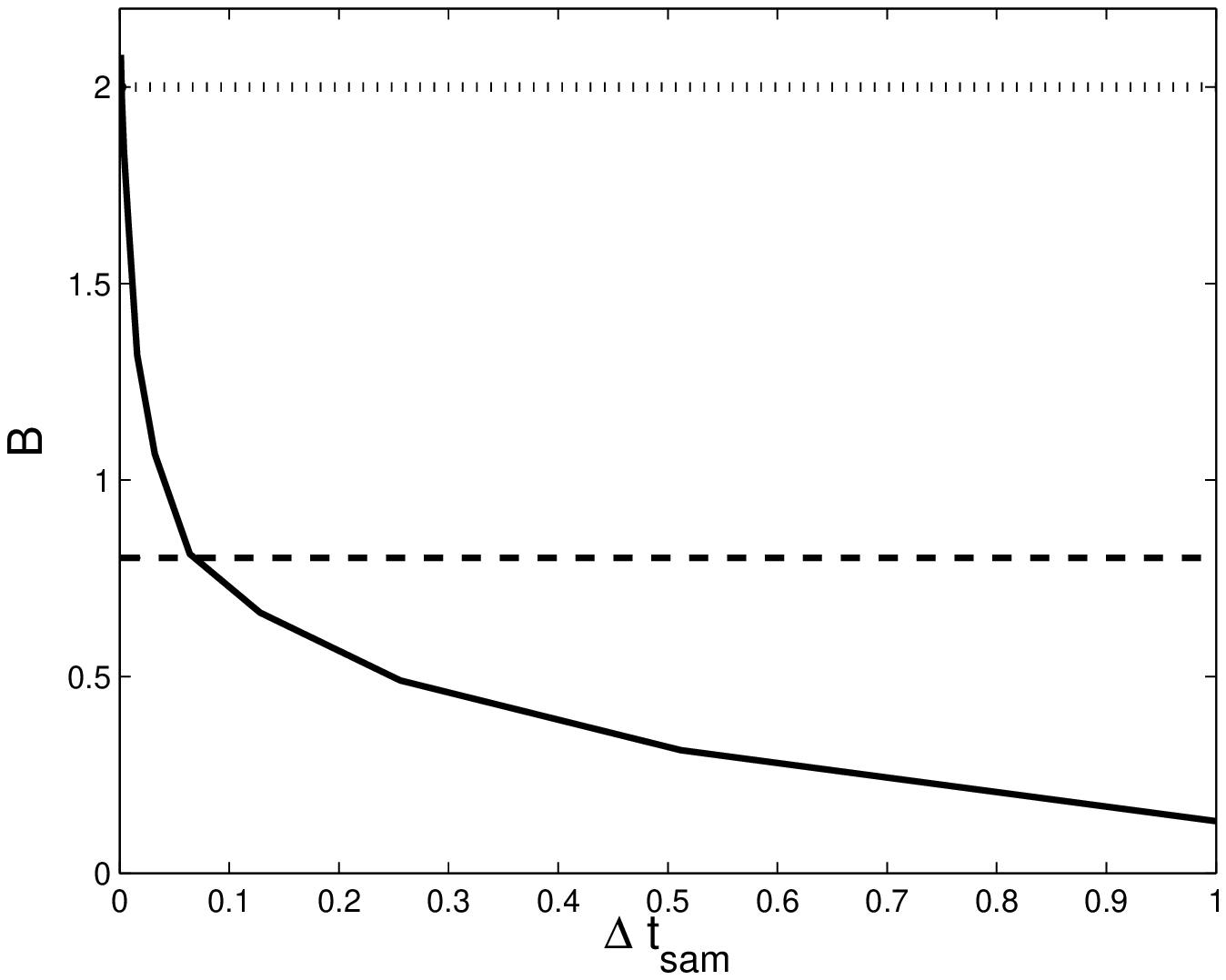} \\
a.~~  $\widehat{A}$ vs $\Delta t_{sam}$  & b.~~ $ \widehat{B}$ vs $\Delta t_{sam}$
\end{tabular}}
\begin{center}
\caption{Estimation of the parameters of the bistable potential \eqref{e:pot_bistable} as a
function of the sampling rate for $\sigma = 0.7, \,\eps = 0.1$. Solid line: estimated coefficient.
Dashed line: homogenized coefficient. Dotted line: unhomogenized coefficient.}
\label{fig:bistable:A_B_07}
\end{center}
\end{figure}

In Figures \ref{fig:bistable:A_B_05} and \ref{fig:bistable:A_B_07} we present the
estimators for the two drift coefficients versus the sampling rate, for two different
values of $\sigma$. We observe that the performance of the estimators is qualitatively
similar to the OU case. Notice also that the optimal sampling rate is
approximately the same for both coefficients.

\begin{figure}[t]
\centerline{
\begin{tabular}{c@{\hspace{2pc}}c}
\includegraphics[width=2.7in, height = 2.7in]
{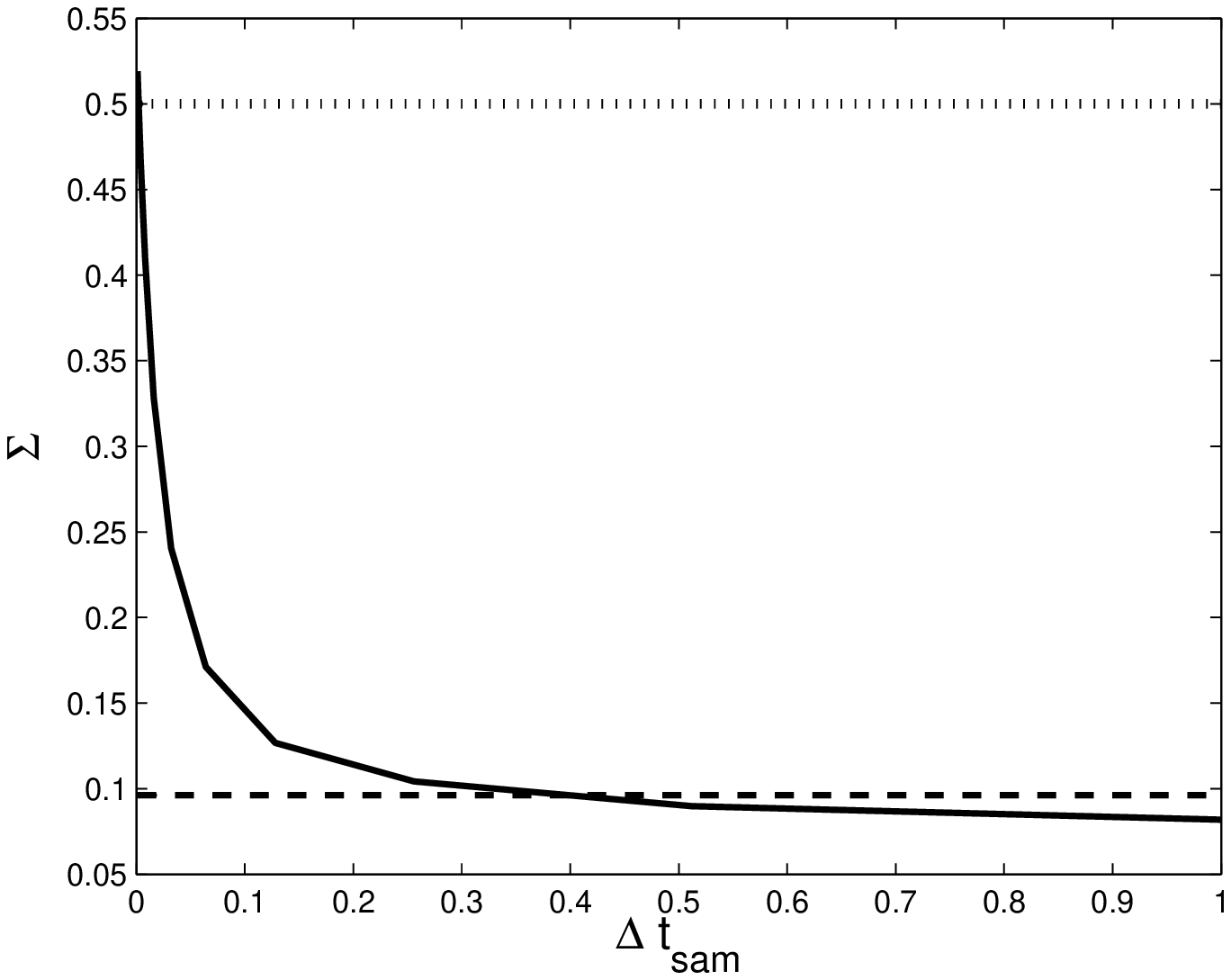} &
\includegraphics[width=2.7in, height = 2.7in]
{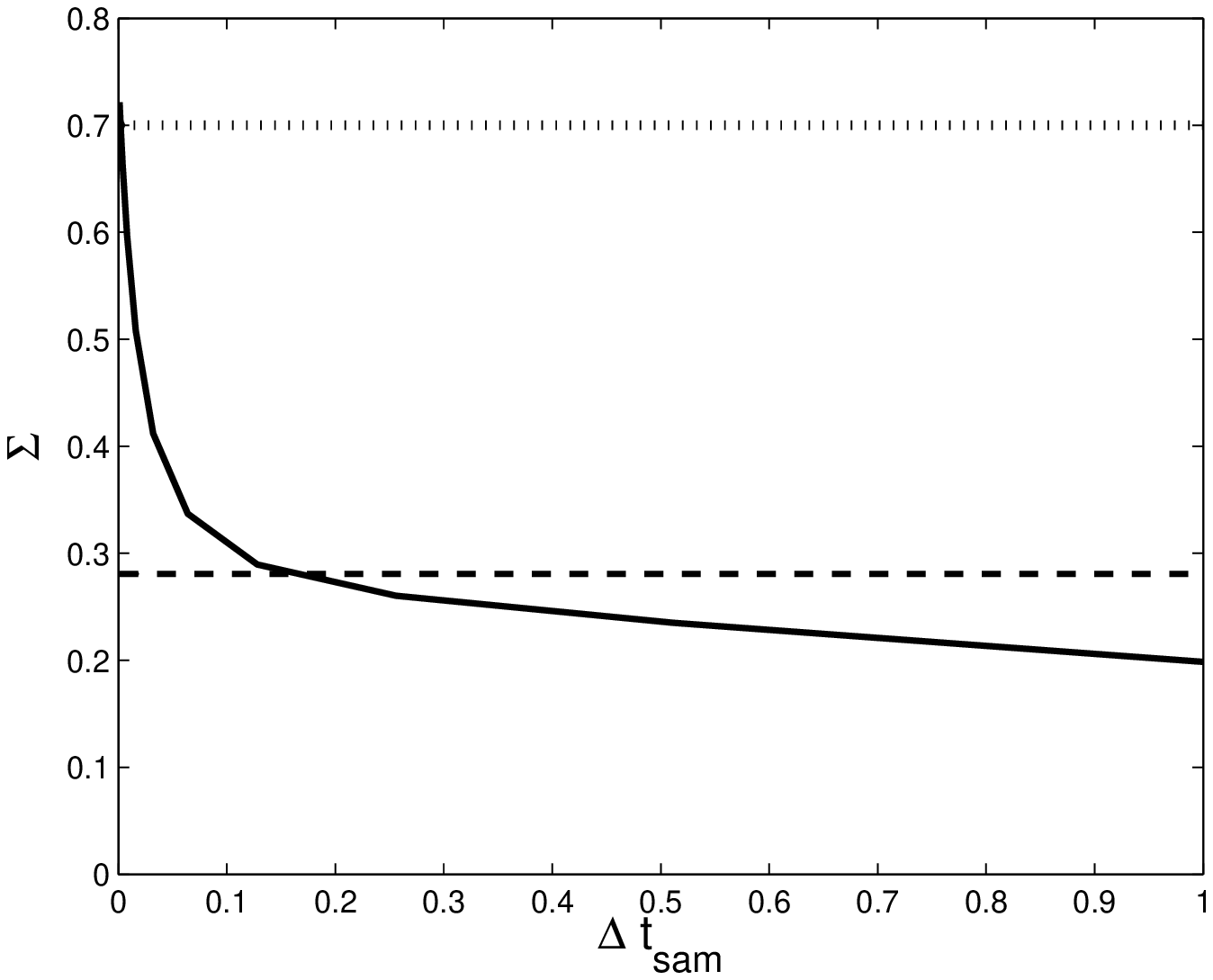} \\
a.~~  $\sigma = 0.5$  & b.~~ $ \sigma =0.7$
\end{tabular}}
\begin{center}
\caption{Estimation of the diffusion coefficient for the bistable potential
\eqref{e:pot_bistable} as a function of the sampling rate for $\alpha = 1.0,
\, \beta = 2.0, \,\eps = 0.1$. Solid line: estimated coefficient. Dashed line:
homogenized coefficient. Dotted line: unhomogenized coefficient.} \label{fig:bistable:sigma}
\end{center}
\end{figure}

In Figure \ref{fig:bistable:sigma} we plot the estimator for the diffusion coefficient
versus the sampling rate, for two different values of $\sigma$. The conclusions reached
from the numerical study of $\widehat{\Sigma}$ for the one dimensional OU process carry
almost verbatim to this case.

%
%%%%%%%%%%%%%%%%%%%%%%%%%%%%%%%%%%%%%%%%%%%%%%%%%%%%%%%%%%%%%%%%%%%%%%%%%%%%%%%%%%%%%%%%%%%
%
%
\subsubsection{A Quadratic Potential in 2D}
We Consider now \eqref{e:main} in two dimensions with a separable fast potential $p(y)$:
\begin{equation}\label{e:2dim}
d x^\eps(t) = - \nabla V(x^\eps(t), B) \, dt - \frac{1}{\epsilon}\nabla p_1 \left(
\frac{x^\eps_1(t)}{\eps} \right)  -\frac{1}{\epsilon}\nabla p_2 \left(
\frac{x^\eps_2(t)}{\eps} \right)  \, dt + \sqrt{2 \sigma } \, d \beta (t),
\end{equation}
where $B$ is the set of the drift parameters that we wish to estimate. The homogenized
equation reads
\begin{equation}\label{e:homog_2d}
d X (t) = - K \nabla V(X (t), B) dt + \sqrt{2 \sigma K} \, d \beta (t),
\end{equation}
where
\begin{equation}\label{e:tensor_2d}
K = \left( \begin{array}{cc}
\frac{L^2}{Z_1 \widehat{Z}_1} & 0  \\
0 & \frac{L^2}{Z_2 \widehat{Z}_2}
\end{array} \right)
\end{equation}
and
\begin{eqnarray*}
Z_i = \int_0^L e^{- \frac{p_i(y_i)}{\sigma}} \, dy_i, \quad
 \widehat{Z}_i = \int_0^L e^{\frac{p_i(y_i)}{\sigma}} \, dy_i, \; \; i=1,2.
\end{eqnarray*}
In the above $L$ denotes the period of $p(y)$.

We will consider the case of a general quadratic potential in two dimensions:
\begin{equation}\label{e:pot_2d}
V(x, B) = \frac{1}{2} x^T B x,
\end{equation}
with $B$ symmetric positive-definite. For the fluctuations we will use a simple
two--dimensional extension of the cosine potential
\eqref{e:cos}:
$$
p_1(y_1) = \cos(y_1), \; p_2(y_2) =  \frac{1}{2}\cos(y_2).
$$
Our goal is to estimate the diffusion tensor and the drift coefficients.
We will estimate the diffusion tensor through the quadratic variation:
\begin{equation}
\widehat{\Sigma}_{N,\delta}(x(t)) = \frac{1}{2 N \delta } \sum_{n = 0}^{N-1}
(x_{n+1} -  x_n ) \otimes (x_{n+1} -  x_n ),
\label{e:sigma_estim_dd}
\end{equation}
where $\otimes$ stands for the tensor product.
For simplicity we will assume that the
diffusion tensor in our model is diagonal. This is consistent with the
homogenized diffusion tensor, see eq.  \eqref{e:tensor_2d}.
We will use generalizations of the maximum likelihood estimator
$\widehat{A}$ in order to estimate the parameters of the quadratic potential.

\begin{figure}[t]
\centerline{
\begin{tabular}{c@{\hspace{2pc}}c}
\includegraphics[width=2.7in, height = 2.7in]{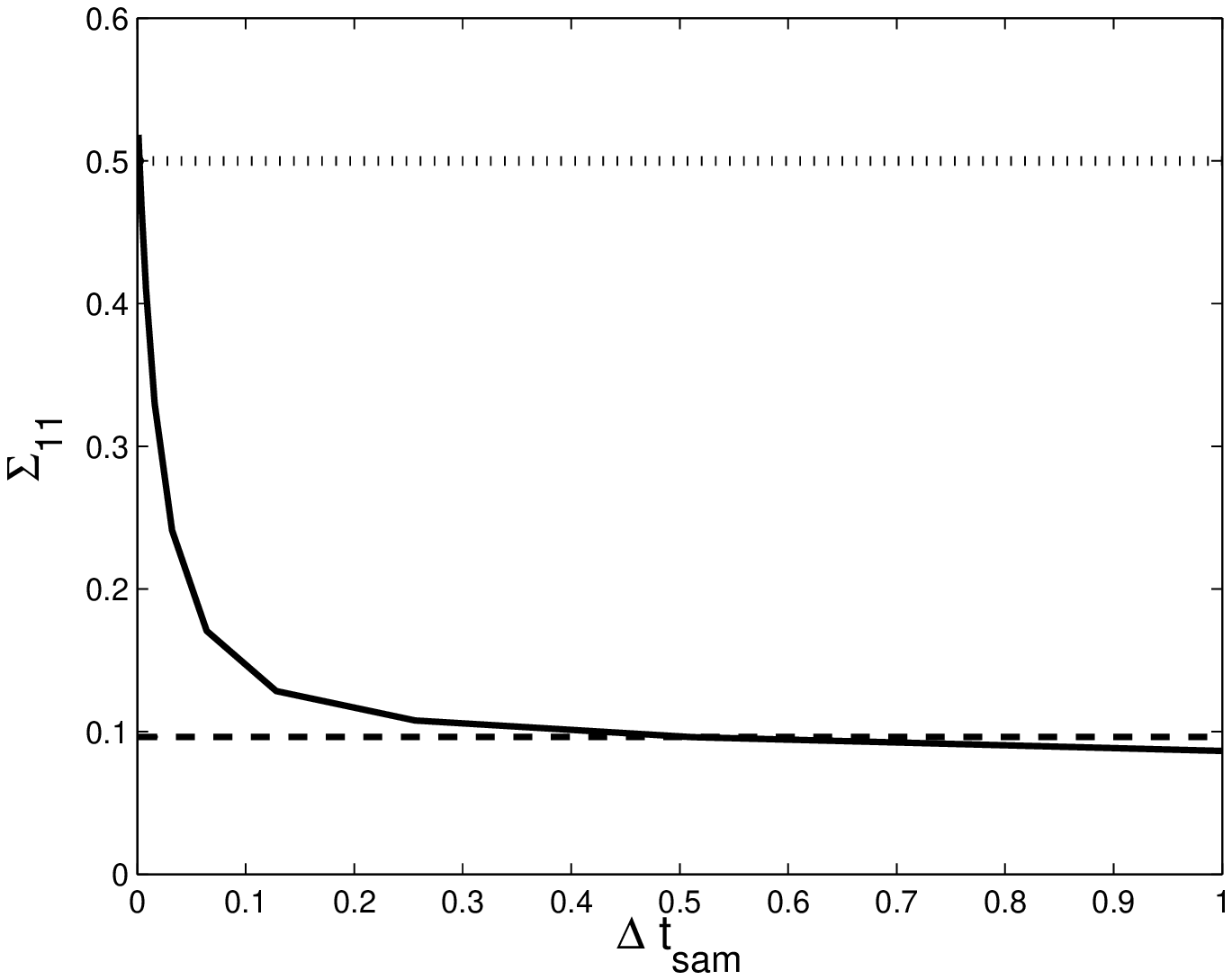} &
\includegraphics[width=2.7in, height = 2.7in]{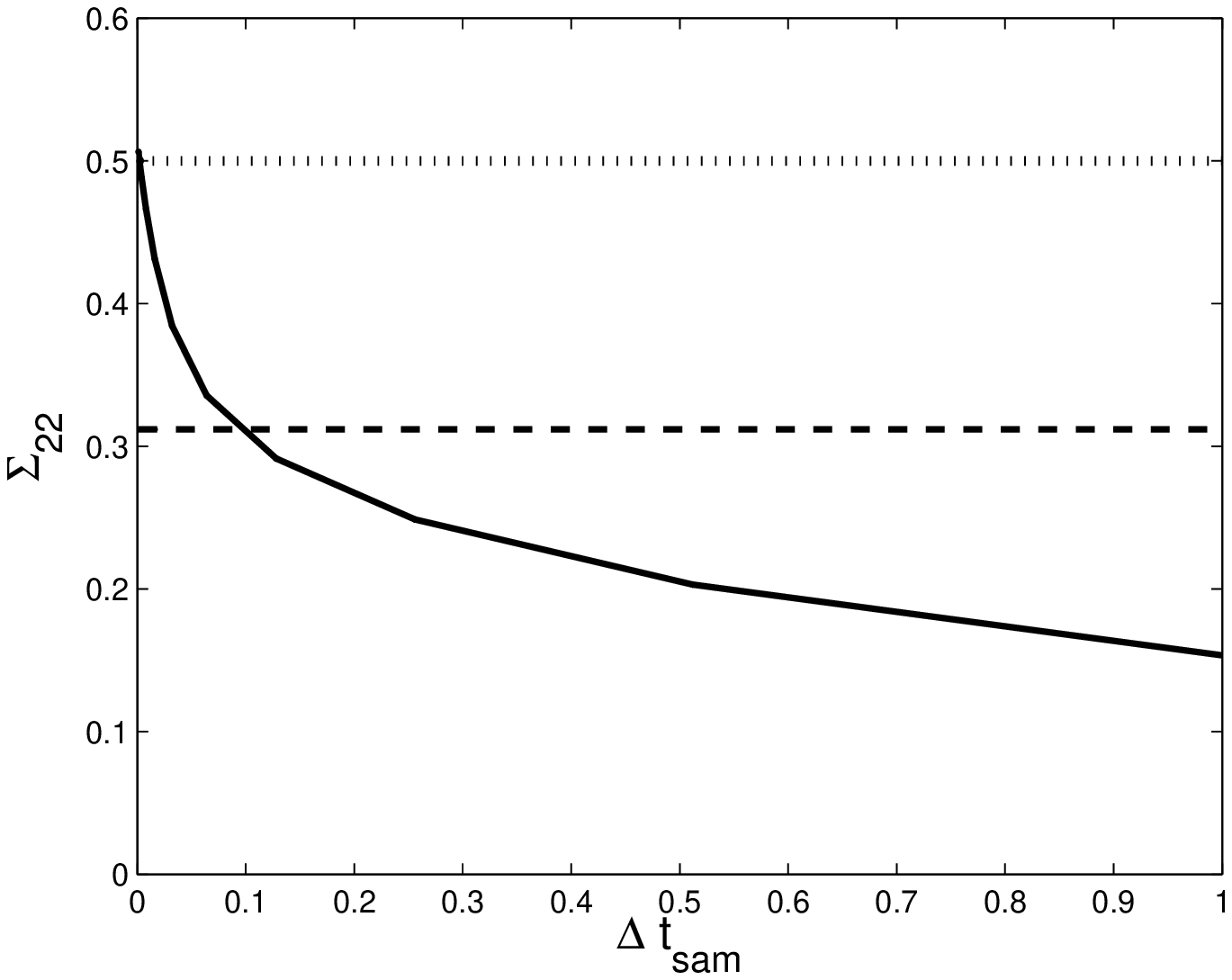} \\
a.~~  $\widehat{\Sigma}_{11}$  & b.~~ $ \widehat{\Sigma}_{22} $
\end{tabular}}
\begin{center}
\caption{Estimation of the non--zero elements of the diffusion tensor for the 2d quadratic potential
\eqref{e:pot_2d} as a function of the sampling rate for $B_{11}= B_{12} = B_{21} = 2, \, B_{22} = 3,
\, \sigma = 0.5, \, \eps = 0.1$. Solid line: estimated coefficient.
Dashed line: homogenized coefficient. Dotted line: unhomogenized coefficient.}
\label{fig:2d_sigma}
\end{center}
\end{figure}
\begin{figure}
\centerline{
\begin{tabular}{c@{\hspace{2pc}}cc}
\includegraphics[width=2.7in, height = 2.7in]{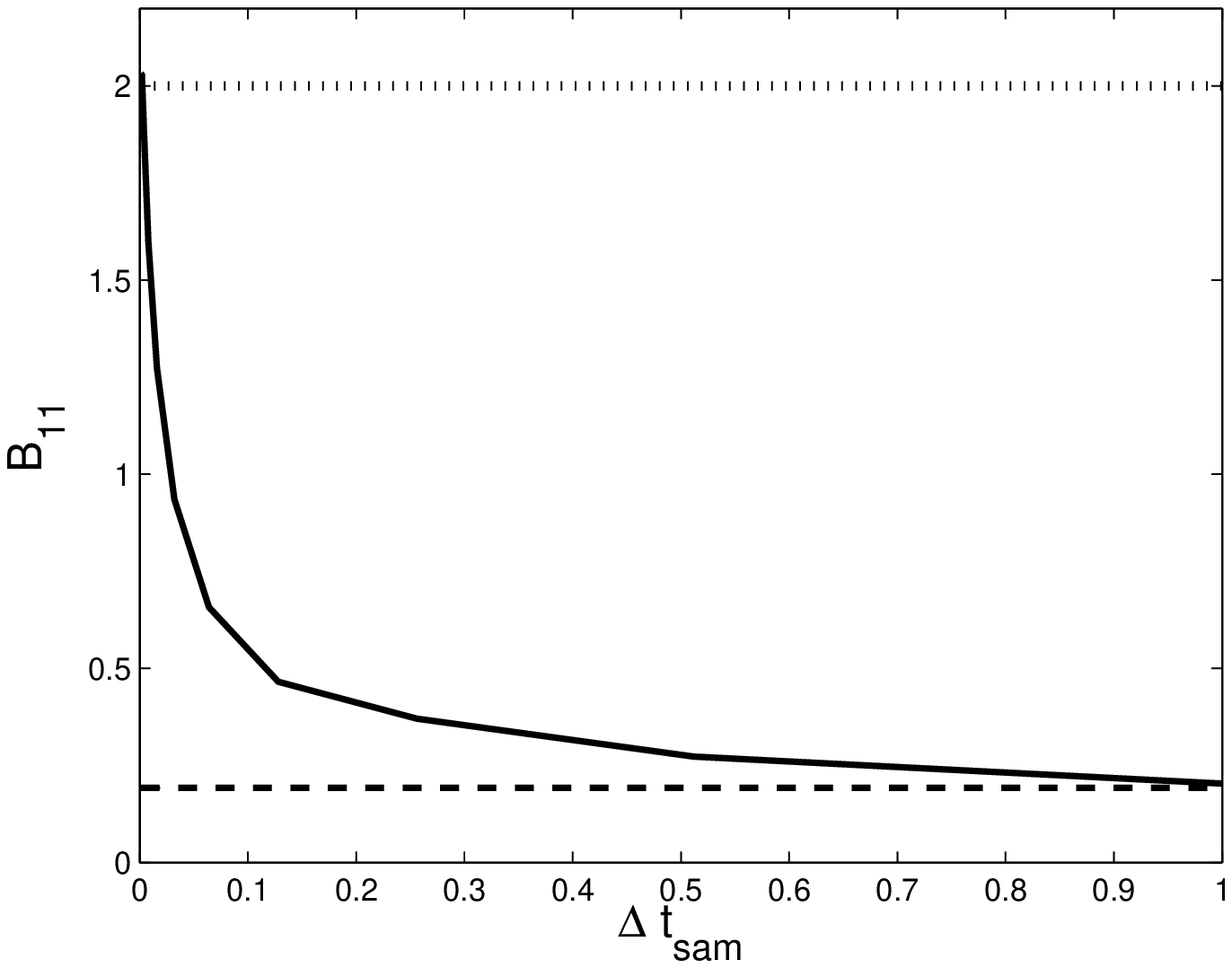} &
\includegraphics[width=2.7in, height = 2.7in]{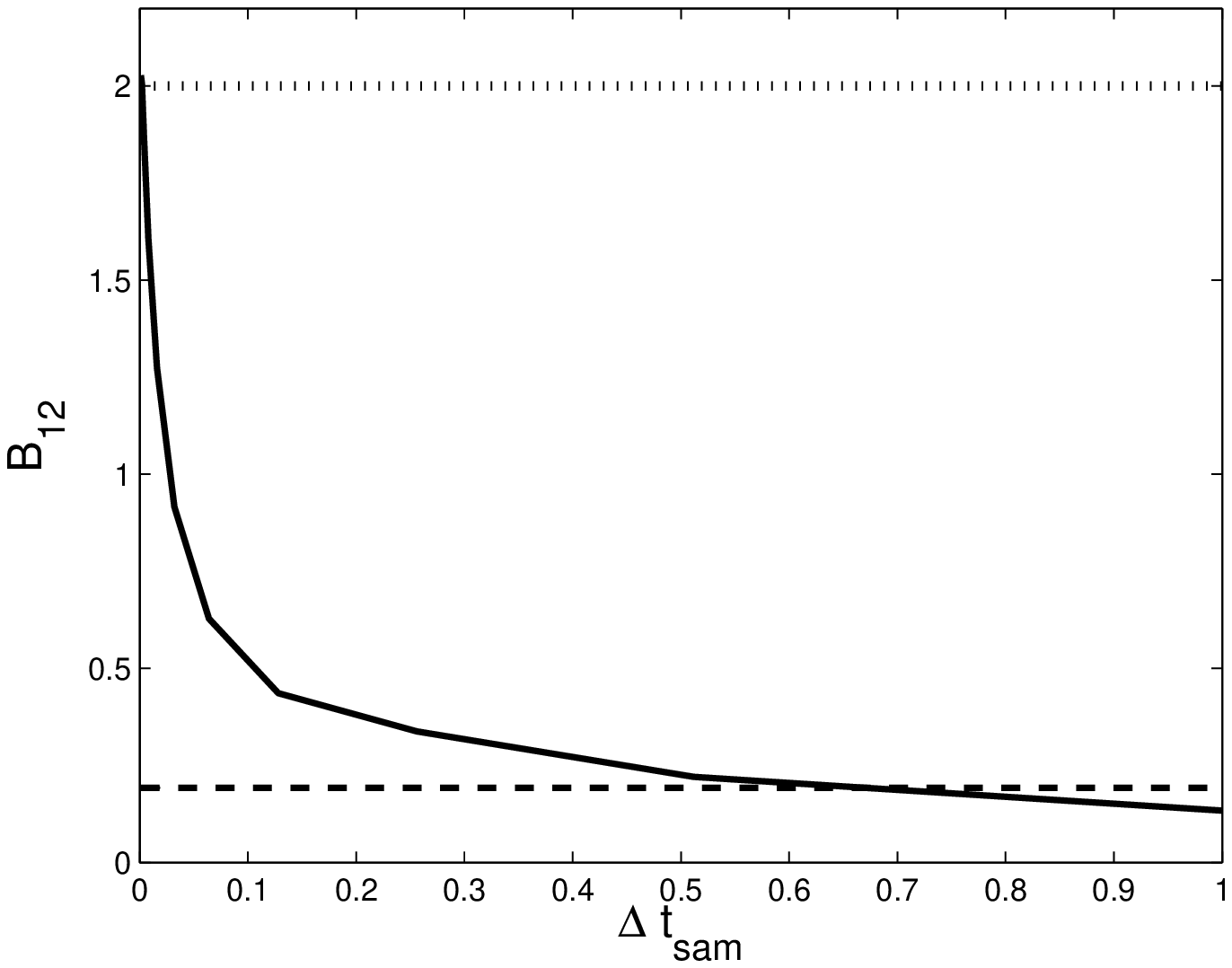} \\
a.~~  $\widehat{B}_{11}$  & b.~~ $ \widehat{B}_{12} $ \\
\includegraphics[width=2.7in, height = 2.7in]{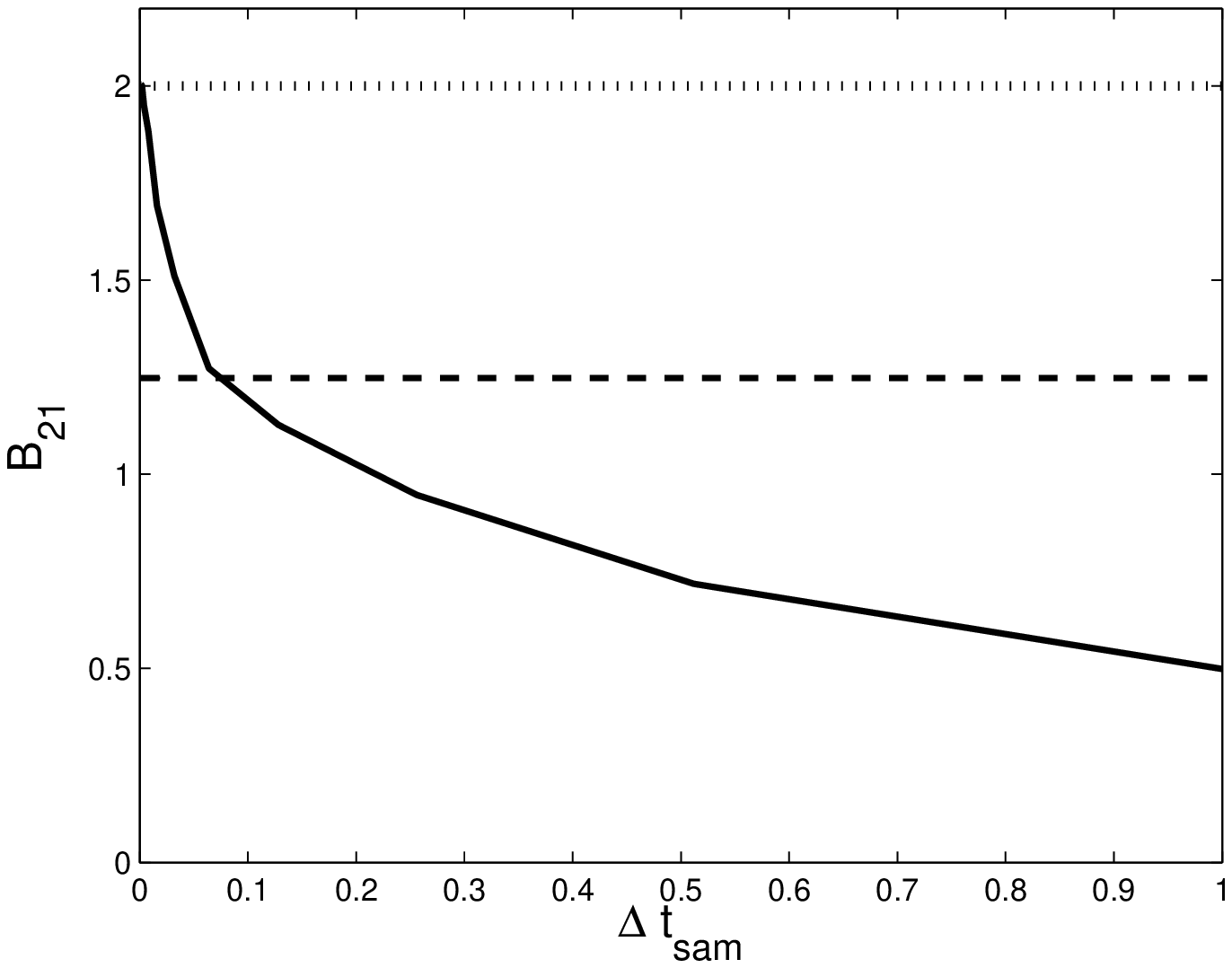} &
\includegraphics[width=2.7in, height = 2.7in]{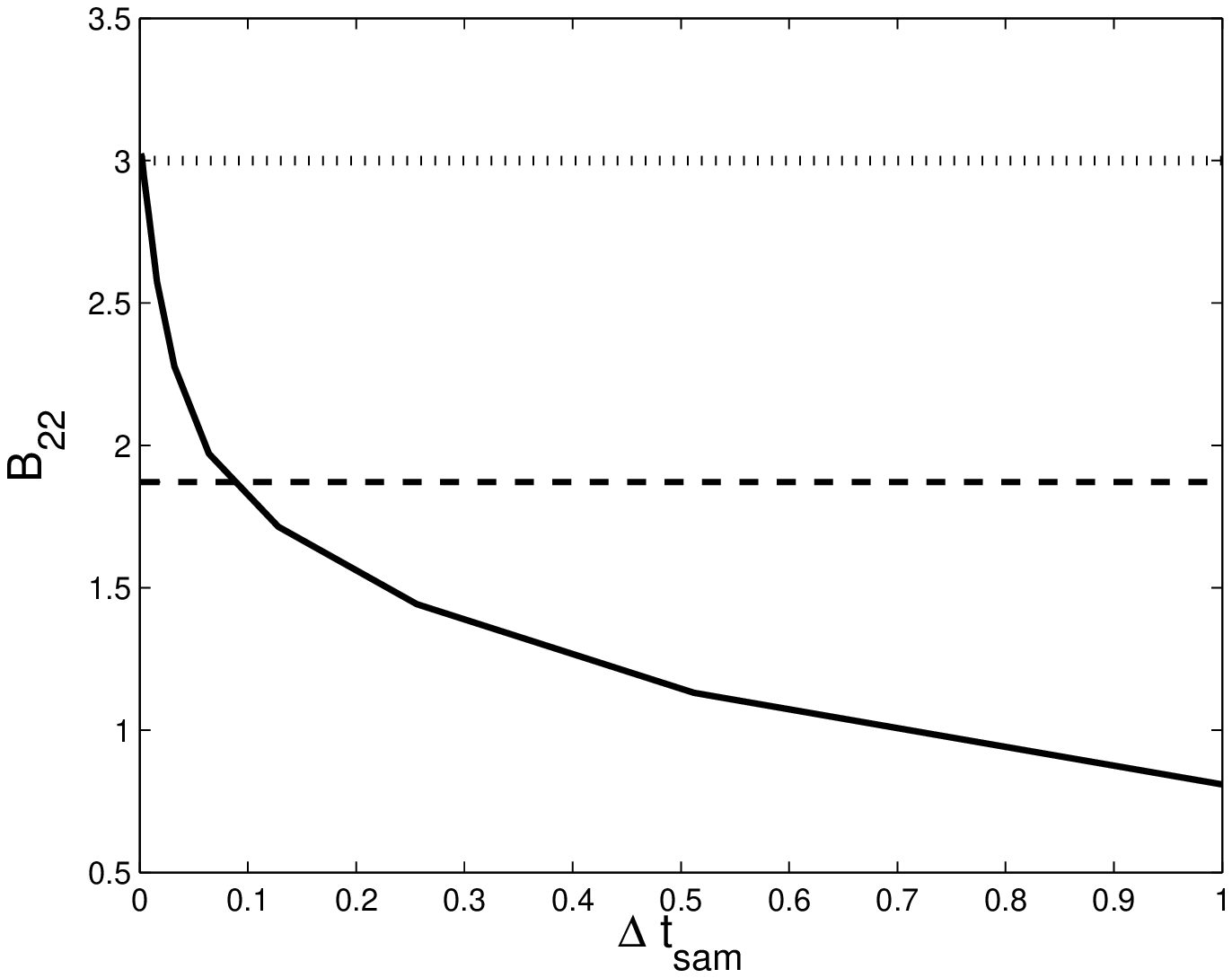} \\
a.~~  $\widehat{B}_{21}$  & b.~~ $ \widehat{B}_{22} $
\end{tabular}}
\begin{center}
\caption{Estimation of the parameters of the 2d quadratic potential
\eqref{e:pot_2d} as a function of the sampling rate for $\sigma = 0.5, \, \eps = 0.1$.
Solid line: estimated coefficient. Dashed line: homogenized coefficient.
Dotted line: unhomogenized coefficient. }
\label{fig:2d_alpha}
\end{center}
\end{figure}
In Figure \ref{fig:2d_sigma} we present the estimated values of the two non--zero
components of the diffusion tensor versus the sampling rate\footnote{The estimated value of
the off--diagonal elements is almost
$0$ for all values of the sampling rate, in accordance with the theoretical result
\eqref{e:tensor_2d}.}. The performance of the estimator for the diffusion tensor is,
qualitatively at least, similar to its performance in the one dimensional problems
considered in the previous two subsections. Notice, however, that the optimal sampling
rate is quite different for the two non--zero components of the diffusion tensor.

In Figure \ref{fig:2d_alpha} we present the estimated values of the four drift
coefficients. The results are in accordance with the one dimensional theory developed in
this paper, as well as with the numerical experiments shown in one dimension. We remark
that the estimators capture successfully the fact that the homogenized matrix $B$ is not
symmetric. Notice furthermore that, as for the diffusion matrix, the optimal sampling
rate is different for different components of the matrix $B$.

Thus, in this simple two dimensional multiscale model, the optimal sampling
rate is different in different directions. This suggests
that extreme care has to be taken when estimating parameters for multidimensional,
multiscale stochastic processes.
%
%
%%%%%%%%%%%%%%%%%%%%%%%%%%%%%%%%%%%%%%%%%%%%%%%%%%%%%%%%%%%%%%%%%%%%%%%%%%%%%%%%%%%%%%%
%
%
\subsection{The Second Estimator for the Drift Coefficient}
\begin{figure}[t]
\centerline{
\begin{tabular}{c@{\hspace{2pc}}c}
\includegraphics[width=2.7in, height = 2.7in]{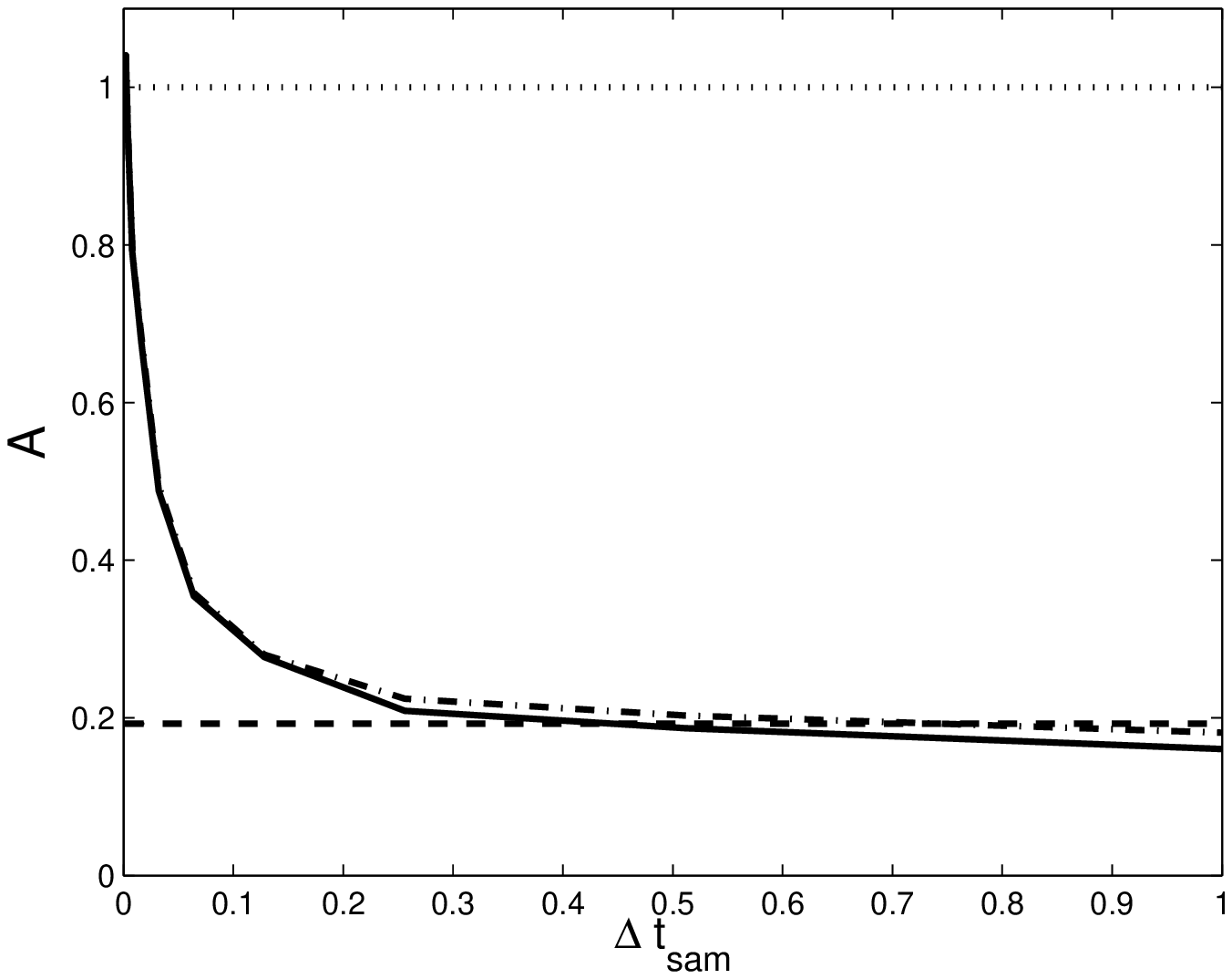} &
\includegraphics[width=2.7in, height = 2.7in]{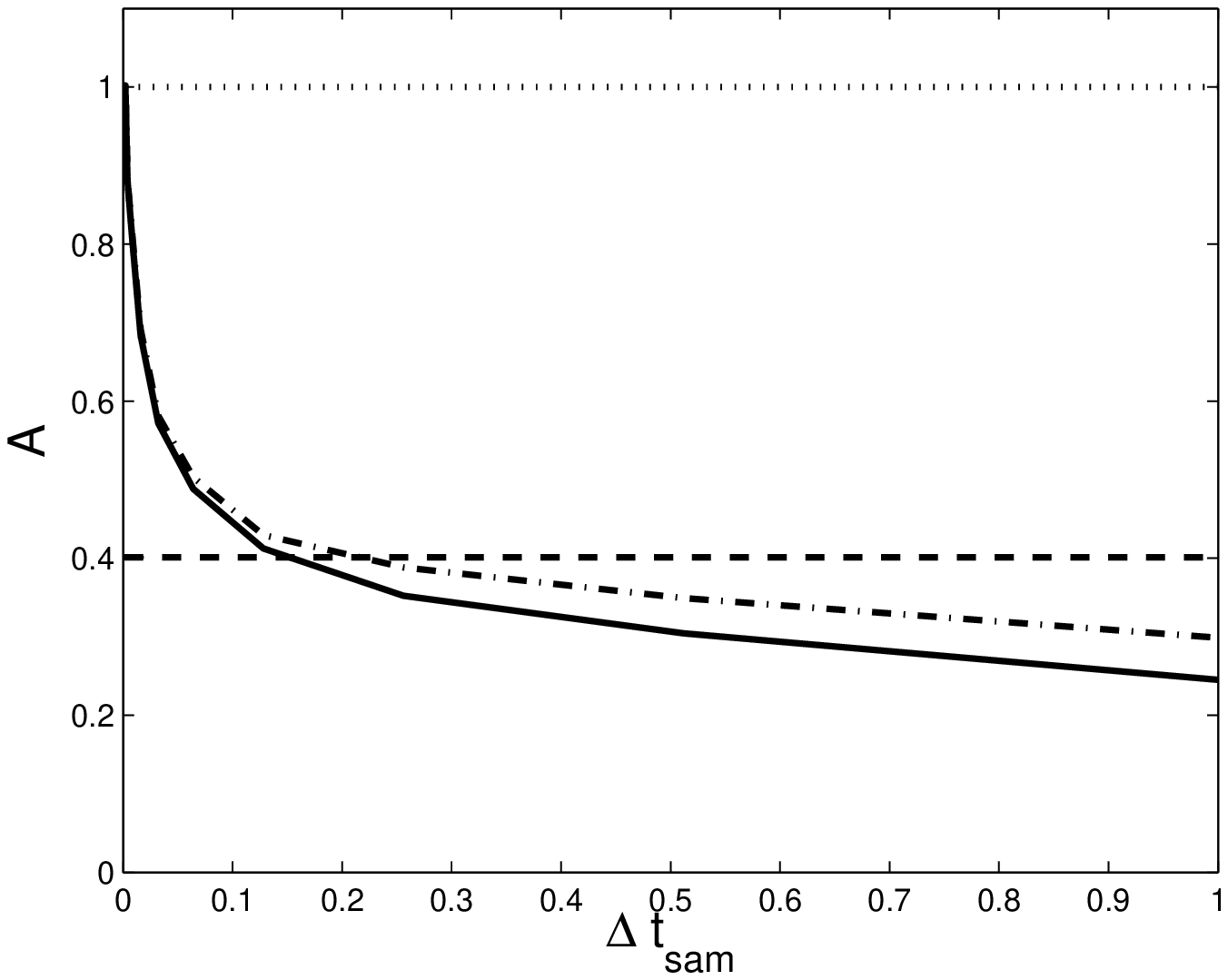} \\
a.~~  $\sigma = 0.5$   & b.~~ $ \sigma = 0.7$
\end{tabular}}
\begin{center}
\caption{ Estimation of the drift coefficients for the quartic potential
\eqref{e:pot_quartic} as a function of the sampling rate for $ \eps = 0.1$. Solid line:
$\widehat{A}$. Dash-dot line: $\tilde{A}$. Dashed line: homogenized coefficient. Dotted
line: unhomogenized coefficient .}
\label{fig:alpha2_four}
\end{center}
\end{figure}

\begin{figure}[t]
\centerline{
\begin{tabular}{c@{\hspace{2pc}}c}
\includegraphics[width=2.7in, height = 2.7in]{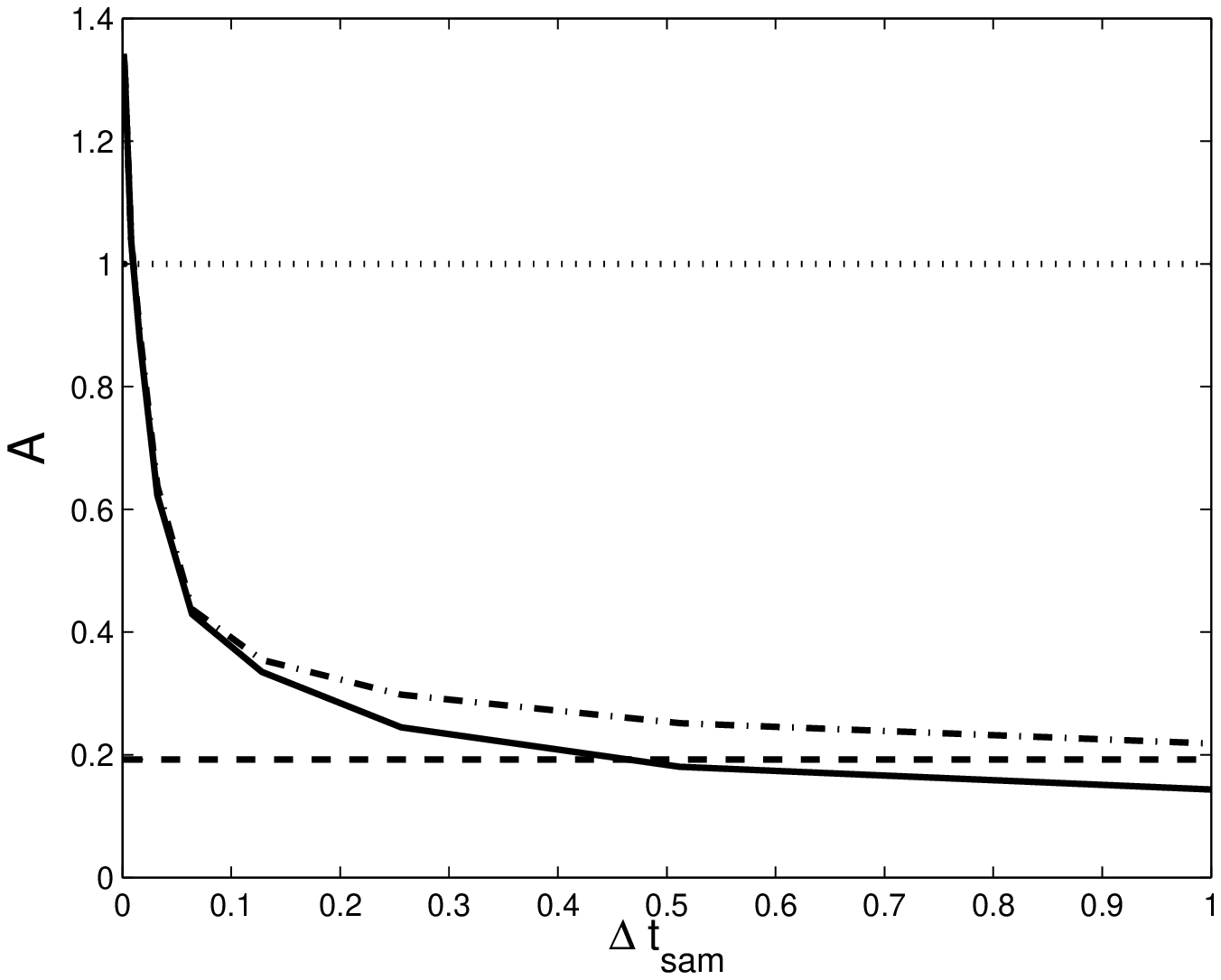} &
\includegraphics[width=2.7in, height = 2.7in]{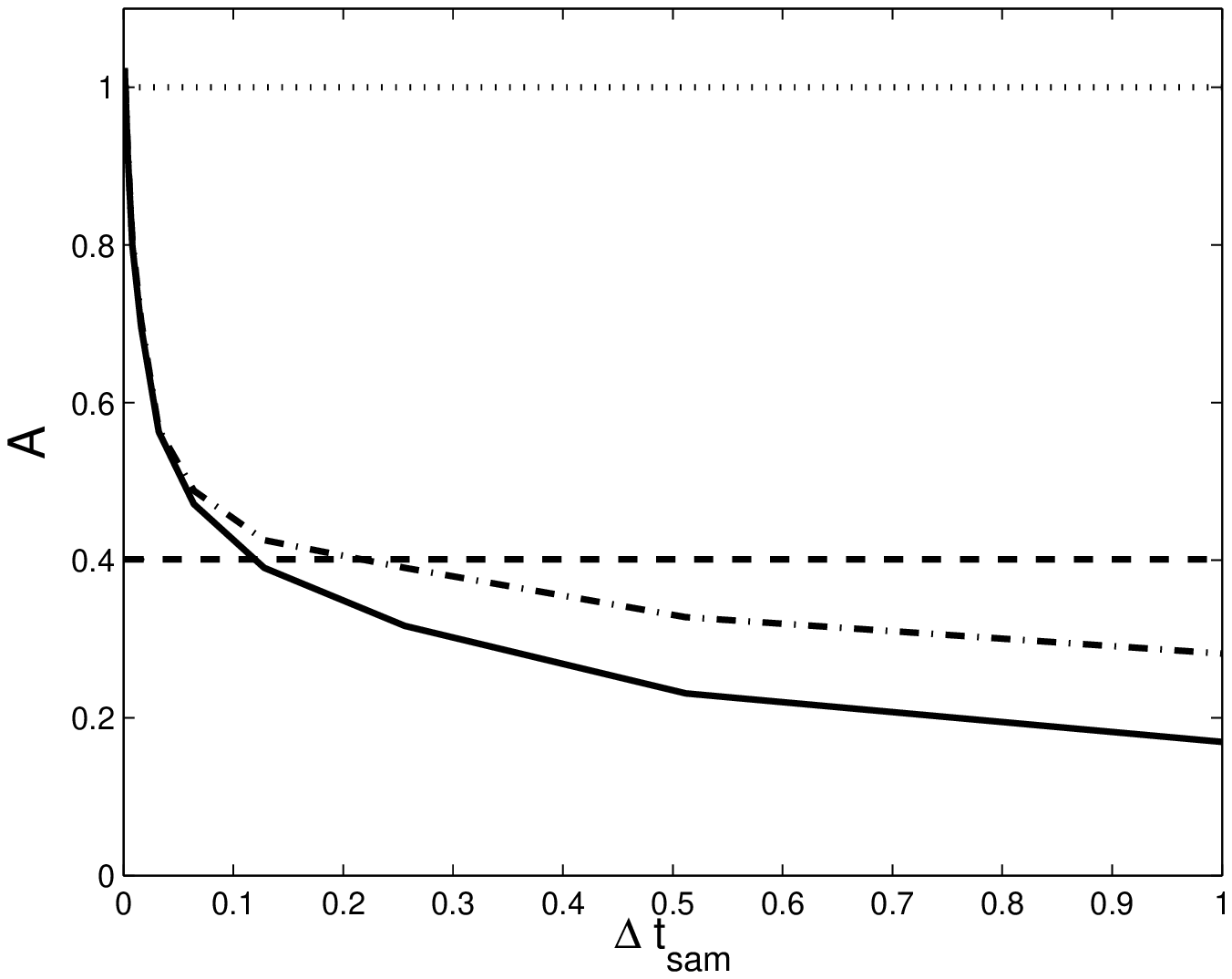} \\
a.~~  $\sigma = 0.5$   & b.~~ $ \sigma = 0.7$
\end{tabular}}
\begin{center}
\caption{ Estimation of the drift coefficients for the sixth--degree potential
\eqref{e:pot_six} as a function of the sampling rate for $\eps = 0.1$. Solid line:
$\widehat{A}$. Dash-dotted line: $\tilde{A}$. Dashed line: homogenized coefficient.
Dotted line: unhomogenized coefficient .}
\label{fig:alpha2_six}
\end{center}
\end{figure}

In this section we compare between the performances of the two estimators
for the drift coefficient, namely $\widehat{A}$ and $\tilde{A}$ given by
equations \eqref{e:alpha_estim_1d} and \eqref{e:alpha_estim_1d2} respectively. We estimate
the drift parameter of
\eqref{e:main} in one dimension for a quartic and a sixth--degree large--scale potential
$V(x)$:
\begin{equation}\label{e:pot_quartic}
V(x) = \frac{1}{4} \alpha x^4
\end{equation}
and
\begin{equation}\label{e:pot_six}
V(x) = \frac{1}{6} \alpha x^6.
\end{equation}
In both cases the small scale fluctuations are represented by the cosine potential
\eqref{e:cos} In Figure \ref{fig:alpha2_four} we present the estimated values of the
drift coefficient as a function of the sampling rate for two different $\sigma$ for the
quartic potential \eqref{e:pot_quartic}. We also plot the effective and the unhomogenized
values of the drift coefficient. Similar results for the sixth--degree potential
\eqref{e:pot_six} are presented in Figure \ref{fig:alpha2_six}. In both cases we observe
that the alternative estimator $\tilde{A}$ performs better than $\widehat{A}$ in this
situation where the data is subsampled.
%
%%%%%%%%%%%%%%%%%%%%%%%%%%%%%%%%%%%%%%%%%%%%%%%%%%%%%%%%%%%%%%%%%%%%%%%%%%%%%%%%
%
%                          STATEMENT OF RESULTS
%
%%%%%%%%%%%%%%%%%%%%%%%%%%%%%%%%%%%%%%%%%%%%%%%%%%%%%%%%%%%%%%%%%%%%%%%%
%
\section{Statement of Main Results}
\label{sec:results}

In this section we pesent theorems which substantiate the numerical
observations in the preceeding section.
The first result shows that, without subsampling, the parameter estimators for the
homogenized model will be asymptotically biased: they recover the parameters from the
unhomogenized equations.
\begin{theorem}
\label{thm:est_ddim} Let $x^\eps(t)$ be the solution of \eqref{e:xeps_V} with $x^\eps
(0)$ distributed according to the invariant measure of the process. Then the estimator
\eqref{e:a_est} satisfies
\begin{equation}\label{e:a_est_lim}
\lim_{\eps \rightarrow 0}\lim_{T \rightarrow \infty} \widehat{A}(x^{\eps}) = \alpha
\quad \mbox{a.s.}
\end{equation}
Fix $T = N \delta$ in \eqref{e:sigma_estim_1d}. Then for every $\eps > 0$ we have
\begin{equation}\label{e:sigma_est_lim}
\lim_{N \rightarrow \infty} \widehat{\Sigma}_{N, \delta}(x^{\eps}) = \sigma \quad
\mbox{a.s.}
\end{equation}
\end{theorem}
Now consider the one dimensional problem
\begin{equation}
d x^\eps (t) = - \alpha V'(x^\eps(t)) dt - \frac{1}{\eps} p' \left(
\frac{x^\eps(t)}{\eps} \right) dt + \sqrt{2 \sigma} d \beta (t). \label{e:xeps_est}
\end{equation}

The next two results show that, with appropriate subsampling, the estimators
recover the correct drift and diffusion coefficients for the homogenized
model \eqref{e:lim_sde_1d} when taking data from the unhomogenized
equation \eqref{e:xeps_est}.
\begin{theorem}\label{thm:par_est_alpha}
Let $x^\eps(t)$ be the solution of \eqref{e:xeps_est} with $x^\eps (0)$ distributed
according to the invariant measure of the process. Further, let
$\delta = \eps^\alpha, \, \alpha \in (0 , 1 )$ and $N = \left[ \eps^{-\gamma} \right], \,
\gamma > \alpha,$ where $[\cdot]$ denotes the integer part of a number. Then
\begin{equation}
\lim_{\eps \rightarrow 0} \widehat{A}_{N, \delta} (x^\eps) = A \quad \mbox{in law,}
\label{e:alpha_lim}
\end{equation}
where $A$ is given by \eqref{e:coeffs_1d}.
\end{theorem}
\begin{theorem}
\label{thm:par_est_sigma} Let $x^\eps(t)$ be the solution of \eqref{e:xeps_est} with
$x^\eps (0)$ distributed according to the invariant measure of the process. Fix $ T = N
\delta$ with $\delta = \eps^\alpha$ and $\alpha \in (0 , 1)$. Then
\begin{equation}
\lim_{\eps \rightarrow 0} \widehat{\Sigma}_{N, \delta} (x^\eps) = \Sigma \quad
\mbox{in law,}
\label{e:sigma_lim}
\end{equation}
where $\Sigma$ is given by \eqref{e:coeffs_1d}.
\end{theorem}
\begin{remark}
The two previous results require $\epsilon/\delta \to 0$ as $\epsilon \to 0.$ In view of
the fact that the fast time--scale is ${\cal O}(\epsilon^2)$ (see equation
\eqref{e:yeps_eqn}) we might expect that this could relaxed to
$\epsilon^2/\delta \to 0$
as $\epsilon \to 0.$ However we have not been able to prove this.
See Remark \ref{r:label} for further discussion of this point.
\end{remark}
The final result concerns the second drift estimator and again concerns
input of data from the unhomogenized equation \eqref{e:xeps_est} into the
paramter estimator for the homogenized equation \eqref{e:lim_sde_1d}.
It requires an estimate of the
diffusion coefficient, $\widehat{\Sigma}.$ If $\widehat{\Sigma} = \sigma$, then we
estimate the drift coefficient incorrectly with $\tilde A(x^{\eps})$; on the other
hand, if $\widehat{\Sigma} = \Sigma$, then the estimator $\tilde A(x^{\eps})$ gives
the drift of the homogenized equation. (To see the last result recall that
$A/\Sigma=\alpha/\sigma$, see \eqref{e:coeffs_1d}). Consequently, for
multiscale gradient systems, it is sufficient only to subsample in a
fashion which leads to the
correct diffusion coefficient. This offers a clear computational advantage.
\begin{theorem}\label{prop:drift_estim_2}
Let $x^\eps(t)$ be the solution of \eqref{e:xeps_est} with $x^\eps (0)$ distributed
according to the invariant measure of the process. Assume that the diffusion coefficient
has been estimated to be $\widehat{\Sigma}$. Then
$$\lim_{\eps \to 0}\lim_{T \to \infty}\tilde A(x^\eps)
=\frac{\widehat\Sigma}{\sigma}\alpha \quad \mbox{in law.}$$
\end{theorem}
%
%
%%%%%%%%%%%%%%%%%%%%%%%%%%%%%%%%%%%%%%%%%%%%%%%%%%%%%%%%%%%%%%%%%%%%%%%%%%%%%%%%%%%%%%%%%%
%
%                                 PRELIMINARY RESULTS
%
%%%%%%%%%%%%%%%%%%%%%%%%%%%%%%%%%%%%%%%%%%%%%%%%%%%%%%%%%%%%%%%%%%%%%%%%%%%%%%%%%%%%%%%%%%%
%
\section{Preliminary Results}
\label{sec:prelim}
In this section we collect various results that will be used in the proof of our main
theorems.
We start by investigating some of the properties of the invariant measures of the
unhomogenized and of the homogenized equation.
We then introduce some tools useful in the study of homogenization for
SDEs.
\begin{prop}
\label{prop:gibbs}
The invariant measure of the homogenized equation \eqref{e:lim_sde}
is the Gibbs measure
\begin{equation}
\mu(dx) = \rho(x) dx = \frac{1}{Z} e^{-\alpha V(x)/\sigma} \, dx,
\quad Z = \int_{\R^d} e^{-\alpha V(x)/\sigma} \, dx.
\label{e:gibbs}
\end{equation}
The Markov process $x(t)$ given by \eqref{e:lim_sde} is geometrically ergodic: there are
$C,\, \lambda>0$ such that, for every measurable $f(x)$ satisfying
$$
|f(x)| \leq 1 + |x|^p,
$$
for some integer $p > 0$, we have, for $\mu-$ a.e. $X(0)$,
$$
\left| \E f(x(t)) - \int_{\R^d} f(x) \rho(x) \, dx  \right| \leq
C\bigl(1+|x(0)^p| \bigr)e^{- \lambda t},
$$
where $\E$ denotes expectation with respect to Wiener measure.
\end{prop}

\proof Assumptions \ref{a:1}, together with the formulae for the effective drift and the
effective diffusion coefficient, equation \eqref{e:coeffs}, imply
that the solution $x(t)$ of the homogenized equation \eqref{e:lim_sde} has a unique
invariant measure with smooth density.  The Gibbs measure \eqref{e:gibbs} satisfies
$$\alpha \nabla V \rho+\sigma \nabla \rho=0$$
and hence
$$K\Bigl(\alpha \nabla V \rho+\sigma \nabla \rho\Bigr)=0.$$
Because $K$ is constant we deduce that
$$\alpha K\nabla V \rho+\nabla \cdot \bigl(\sigma K \rho\bigr)=0.$$
Thus
$$\nabla \cdot \Bigl(\alpha K\nabla V \rho+\nabla \cdot \bigl(\sigma K \rho\bigr)\Bigr)=0.$$
This is the stationary Fokker-Planck equation for \eqref{e:lim_sde}
showing that the Gibbs measure $\rho$ is indeed an invariant measure.
For the geometric ergodicity we use \cite[Thm 5.3]{MattStuHigh02}.
\qed
\begin{prop}\label{lem:xeps_meas_ddim}
The invariant measure of the unhomogenized equation \eqref{e:xeps_V} is the Gibbs measure
\begin{equation}
\mu^\eps(dx) = \rho^\eps(x) \, dx = \frac{1}{Z^\eps} e^{-\frac{\alpha}{ \sigma}V(x) -
\frac{1}{\sigma} p \left(\frac{x}{\eps} \right)}, \quad Z^\eps := \int_{\R^d}
e^{-\frac{\alpha}{ \sigma} V(x) - \frac{1}{\sigma} p \left(\frac{x}{\eps} \right)} \, dx.
\label{e:xeps_inv_meas_ddim}
\end{equation}
For every $\eps > 0$ the Markov process \eqref{e:xeps_V} is geometrically
ergodic: there are $C,\lambda>0$ such that, for every measurable $f(x)$
satisfying
$$
|f(x)| \leq 1 + |x|^p,
$$
for some integer $p>0$ we have, for $\mu^{\eps}-$a.e. $x^{\eps}(0)$,
$$
\left| \E f(x^\eps(t)) - \int_{\R} f(x) \rho^\eps(x) \, dx  \right|
 \leq C\bigl(1+|x^{\eps}(0)|^p\bigr)e^{- \lambda t},$$
where $\E$ denotes expectation with respect to Wiener measure.

Furthermore, the measure $\mu^\eps$ converges weakly to the invariant measure of the
homogenized dynamics $\mu$ given by \eqref{e:gibbs}.
\end{prop}
\proof Assumptions \ref{a:1} imply that $x^\eps(t)$ is an ergodic Markov process. Direct
calculation with the Fokker--Planck equation shows that the unique invariant measure of
the process is the Gibbs measure
\begin{eqnarray*}
\rho^\eps(x) \, dx & = & \frac{1}{Z^\eps} e^{- \frac{1}{\sigma}
 V \left( x, \frac{x}{\eps}, \alpha \right)} \, dx
\\ & = & \frac{1}{Z^\eps} e^{- \frac{\alpha}{\sigma}
 V(x)- \frac{1}{\sigma} p\left(  \frac{x}{\eps} \right)} \, dx,
\end{eqnarray*}
with $Z^\eps$ given by \eqref{e:xeps_inv_meas_ddim}.
For the geometric ergodicity we use \cite[Thm 5.3]{MattStuHigh02}.

Now let
$$
u(x,y):= e^{- \frac{\alpha}{ \sigma} V(x) - \frac{1}{\sigma}p(y)}.
$$
Since $u(x,y) \in L^1(\R^d ; C_{per}(\T^d))$, by \cite[Lem. 9.1]{cioran} we have that
$$
u \left(\cdot, \frac{\cdot}{\eps}  \right) \rightharpoonup \int_{\T^d}
u(\cdot, y) \, dy, \quad \mbox{weakly in } L^1(\R^d).
$$
In particular, since $1 \in L^{\infty}(\R^d)$,
$$
\lim_{\eps \rightarrow 0} Z^\eps =  \int_{\R^d} \int_{\T^d} e^{-
\frac{\alpha}{\sigma}V(x) - \frac{1}{\sigma} p(y)} \, dy.
$$
We combine the above two results to conclude that
$$
\rho^\eps(x) \rightharpoonup \frac{1}{Z} e^{-\frac{\alpha}{ \sigma} V(x)}, \quad
\mbox{weakly in } L^1(\R^d),
$$
where $Z$ is given by \eqref{e:gibbs}. The weak convergence of the densities in
$L^1(\R^d)$ implies the weak convergence of the corresponding probability measures. \qed

\begin{remark}
The assumption of stationarity of the process $x^\eps(t)$ is not necessary for the proof of
the above theorems and is only made for simplicity. Indeed, in the next section we prove that
$x^\eps (t)$ is geometrically ergodic and consequently it converges to its invariant
distribution exponentially fast for arbitrary initial conditions.
Furthermore, the fact that the invariant measure of the process
$x^\eps(t)$ converges weakly, as $\eps \rightarrow 0$,
to the invariant measure of the homogenized process is important for
us as many of our results will be deduced by taking expectations with respect to the
invariant measure $\mu^{\eps}(dx)$ of the multiscale dynamics \eqref{e:xeps_V}. The weak
convergence alluded to demonstrates that the measure $\mu^{\eps}$ behaves uniformly in
$\eps \to 0.$
\end{remark}

An immediate corollary of the above proposition is that $x^\eps(t)$ has bounded moments
of all orders. We will use the notation $\bbE^{\mu^{\epsilon}}$ to denote expectation
with respect to the stationary measure of \eqref{a:1} on path space, when initial data is
distributed according to the Gibbs measure \eqref{e:xeps_inv_meas_ddim}.
\begin{corollary}\label{cor:moments}
Let $x^\eps(t)$ be the solution of \eqref{e:main} with the potential given by
\eqref{e:potential} and assume that conditions \eqref{a:1} are satisfied. Assume
furthermore that $x^\eps(0)$ is distributed according to $\mu^\eps$. Then, for all $p \ge
1,$ there is a constant $C=C(P,T)$ uniform in $\epsilon \to 0$, such that
$$\bbE^{\mu^{\epsilon}}|x^\eps(t)|^p \le C \quad \forall \, t \in [0,T].$$
\end{corollary}

It is convenient for the subsequent analysis to introduce the auxiliary
variable
$$
y^\eps (t) = \frac{x^\eps(t)}{\eps}.
$$
We can then write equation \eqref{e:xeps_V} in the form
\begin{subequations}
\begin{equation}
d x^\eps(t) = - \alpha \nabla V(x^\eps(t)) \, dt -  \frac{1}{\eps} \nabla p \left( y^\eps
(t) \right) \, dt + \sqrt{2 \sigma} \, d \beta (t),
\label{e:xeps_eqn}
\end{equation}
\begin{equation}
d y^{\eps}(t) = - \frac{1}{\eps} \alpha \nabla V(x^\eps(t)) \, dt -  \frac{1}{\eps^2}
\nabla p \left( y^\eps (t)  \right) \, dt + \sqrt{\frac{2 \sigma}{\eps^2}} \, d \beta
(t).
\label{e:yeps_eqn}
\end{equation}
\label{e:eqns_motion}
\end{subequations}
Notice that both processes $x^\eps (t)$ and $y^\eps (t)$are driven by the same Brownian
motion. Written in this fashion it is clear that we are in a situation
where homogenization applies. The homogenized equation is found by
eliminating $y^{\eps}(t)$ from the scale separated system
for $\left\{ x^{\eps}(t), y^{\eps}(t) \right\}$. Note that
${\cal L}_0$ defined in \eqref{e:cell} is the generator of the process
\begin{equation*}
d y (t) = - \nabla p \left( y (t)  \right) \, dt + \sqrt{2 \sigma} \, d \beta (t),
\end{equation*}
on the unit torus, which governs the dynamics of $y_t^{\eps}$ to leading order in
$\epsilon$. The generator of the joint process $\{x^\eps(t), \, y^\eps_t \}$ reads
$$\LL^\eps=\frac{1}{\eps^2} \LL_0 + \frac{1}{\eps} \LL_1 + \LL_2,$$
where
\begin{align*}
\LL_0&= - \nabla_y p(y)\cdot \nabla_y + \sigma \Delta_y,\\
\LL_1&=  - \nabla_y p(y) \cdot \nabla_x - \alpha \nabla_x V(x)
\cdot \nabla_y + 2 \sigma \nabla_x \cdot \nabla_y,\\
\LL_2&=  - \alpha \nabla_x V(x) \cdot \nabla_x + \sigma \Delta_x.
\end{align*}

The following result can be found in, e.g. \cite[Ch. 3]{lions}.
\begin{lemma}
\label{l:Poisson} Assume that $p(y) \in C^{\infty}_{per}(\T^d,\R)$ and that $H(y) \in
C^{\infty}_{per}(\T^d,\R^d).$ Let $\mu(dy)$ be the Gibbs measure \eqref{e:gibbs_torus}
and assume that $H(y)$ is centered with respect to $\mu(dy)$:
\begin{equation}\label{e:centering}
\int_{\T^d} H(y) \, \mu(dy) = 0.
\end{equation}
 Then the Poisson equation
\begin{equation}\label{e:poisson}
- \LL_0 \chi = H(y),
\end{equation}
has a unique mean-zero solution in $L^2_{per}(\T^d, \mu(dy) ; \R^d)$.
This solution, together with all its derivatives, is bounded.
\end{lemma}

We will need an estimate on integrals whose integrand is centered with respect to the
invariant measure  $\mu(dy)$.
\begin{lemma}\label{lem:ito}
Let $H(y) \in C^\infty_{per} \left(\T^d ; \R^d \right)$ satisfy condition
\eqref{e:centering}.  Assume that  $x^\eps (0)$ is distributed according
to \eqref{e:xeps_inv_meas_ddim}. Then the
following estimate holds for any $p>1$ and $T >0$:
\begin{equation*}
 \E^{\mu^\eps} \left| \int_{0}^{T} H(y^\eps(s)) \, ds \right|^p \leq
C \left(\eps^{2p} + \eps^pT^p+\eps^p T^{\frac{p}{2}}  \right).
\end{equation*}
\end{lemma}
\proof Consider the Poisson equation \eqref{e:poisson} with periodic boundary conditions.
Since $H(y)$ satisfies \eqref{e:centering}, Lemma \ref{l:Poisson} applies and we have
that $\chi(y)$ is smooth and bounded, together with all its derivatives. We now apply the
It\^{o} formula to $\chi(y^\eps (t))$, where $y^\eps (t)$ is the solution of
\eqref{e:yeps_eqn}, and use \eqref{e:poisson} to obtain
\begin{align*}
\int_{0}^{T} H(y^\eps(s)) \, ds =& - \eps^2 \left( \chi(y^\eps(T))  -
\chi(y^\eps(0))   \right)\\
& + \eps \sqrt{2 \sigma} \int_{0}^{T} \langle \nabla_y \chi(y^\eps(s)), \, d \beta(s)
\rangle -\alpha \eps \int_0^T \langle \nabla V(x^\eps(s)),
\nabla \chi(y^\eps(s)) \rangle ds.
\end{align*}
Now, using the boundedness of $\chi$, we have, for
$$I(T) :=\E^{\mu^\eps} \left| \int_{0}^{T} H(y^\eps(s)) \, ds, \right|^p,$$
\begin{eqnarray*}
I(T) & \leq & C \left( \eps^{2 p} + \eps^{p} \E^{\mu^\eps}\left|\int_0^T |\nabla
V(x^{\eps}(s))| ds\right|^p+ \eps^{p} \E^{\mu^\eps} \left| \int_{0}^{T} \langle \nabla_y
\chi(y^\eps(s)) , \, d \beta (s) \rangle \right|^p \right)
\\ & \leq &
 C \left( \eps^{2 p} + \eps^{p} T^{p-1} \int_0^T |x^{\eps}(s)|^p ds+
\eps^{p} T^{\frac{p}{2} -1} \int_{0}^{T} \E^{\mu^\eps} \left| \nabla_y \chi(y^\eps(s))
\right|^p \, ds \right)
\\ & \leq &
C \left( \eps^{2 p} + \eps^p T^p+\eps^{p} T^{\frac{p}{2}} \right),
\end{eqnarray*}
from which the desired estimate follows. In deriving the above we used
the estimate \cite[Eqn. 3.25, p. 163]{KSh91} on moments of
stochastic integrals. \qed

%The above result will be of particular use to us in the case $T=\delta \ll 1.$
%The Markov property then implies that
%\begin{equation}\label{e:ito_est}
%\left( \E^{\mu^\eps} \left| \int_{n\delta}^{(n+1)\delta} H(y^\eps(s)) \, ds
%\right|^p \right)^{1/p} \leq   C \left(\eps^{2p} + \eps^p\delta^{\frac{p}{2}}
%\right)^{1/p}. %\end{equation}
%
For the rest of this section we will restrict ourselves to the one dimensional case.
 If we apply It\^{o} formula to $\phi(y^\eps(s))$, the solution of the Poisson equation
\eqref{e:cell}, then we obtain
\begin{eqnarray}
x^\eps_{n+1} - x_n^\eps & = & - \alpha \int_{n \delta}^{(n+1)\delta} V'(x^\eps(s)) (1 +
\partial_y \phi(y^\eps(s))) \, ds
\\ &&+ \sqrt{2 \sigma} \int_{n \delta}^{(n+1)\delta} (1 +
\partial_y \phi(y^\eps(s))) \, d \beta (s)
\nonumber\\ && - \eps \left( \phi(y^\eps((n+1)\delta)) - \phi(y^\eps (n\delta)) \right).
\label{e:integr_parts}
\end{eqnarray}
The proof Theorems \ref{thm:par_est_alpha} and \ref{thm:par_est_sigma} is based on
careful asymptotic analysis of the behavior of $x^\eps_{n+1} - x^\eps_n$ given by this
formula when both $\eps$ and $\delta$ are small. Specifically we will use the following
two propositions. They show how the effective homogenized behaviour is manifest in
the time--$\delta$ Markov chain induced by sampling the path $x^{\eps}(t)$ from
\eqref{e:xeps_V}.
\begin{prop}\label{prop:xndelta1}
For $\eps, \, \delta >0 $ sufficiently small and $n \in \mathbb{N}$ there exists an i.i.d.
sequence of random variables $\xi_n \in \mathcal{N}(0,1)$ such that
\begin{equation}
\sqrt{2 \sigma} \int_{n \delta}^{(n+1) \delta}(1+\partial_y \phi(y^\eps(s))) \, d \beta(s)
=\sqrt{2 \Sigma \, \delta} \, \xi_n+ R_1(\delta, \eps)
\label{e:xn_loc1}
\end{equation}
in law. The remainder $R_1(\delta, \eps)$ satisfies, for every $\beta \in
(0,\frac12)$ and $p>0$, the estimate
\begin{equation}
\left( \E^{\mu^\eps} \big| R_1(\eps,\delta) \big|^p \right)^{1/p} \leq  C \, \left(
\eps^{2 \beta} + \eps^{\beta} \right),
\label{e:R_est}
\end{equation}
where $C$ is independent of $\eps$ and $\delta$.
\end{prop}
\begin{remark}
\label{r:label}
Estimate \eqref{e:R_est} is almost certainly not optimal. Indeed, informal
calculations lead us to expect the estimate
$$
\left( \E^{\mu^\eps} \big| R_1(\eps,\delta) \big|^p \right)^{1/p} \leq  C \, \left(
\eps^{2 \beta} + \eps^{\beta} \delta^{\beta} + \eps^{\beta} \delta^{\frac{\beta}{2}} \right).
$$
However, we have not been able to prove this.
\end{remark}
\begin{prop}\label{prop:xndelta2}
For $\eps, \, \delta >0 $ sufficiently small and $n \in \mathbb{N}$ we have that
\begin{equation}
\label{e:xn_loc2}
\alpha \int_{n \delta}^{(n+1) \delta} V'(x^\eps (s)) (1 +
\partial_y \phi(y^\eps (s))) \, ds = A \delta V'(x^\eps_n) +R_2(\eps, \delta)
\end{equation}
in law. The remainder $R_2(\delta, \eps)$ satisfies, for every $p>0$,
the estimate
\begin{equation}
\left( \E^{\mu^\eps} \big| R_2(\eps,\delta) \big|^p \right)^{1/p} \leq  C \left(
\eps^2 +\delta^{ \frac12}\eps+ \delta^{3/2} \right), \label{e:R_est2}
\end{equation}
where $C$ independent of $\eps$ and $\delta.$
\end{prop}

%%%%%%%%%%%%%%%%%%%%%%%%%%%%%%%%%%%%%%%%%%%%%%%%%%%%%%%%%%%%%%%%%%%%%%%%%%%%%%%%%%%%%%%%%%%%
%
%                   PROOF OF PROPOSITION 1.1
%
%%%%%%%%%%%%%%%%%%%%%%%%%%%%%%%%%%%%%%%%%%%%%%%%%%%%%%%%%%%%%%%%%%%%%%%%%%%%%%%%%%%%%%%%%%%
%
\section{Proof of Propositions \ref{prop:xndelta1} and \ref{prop:xndelta2}}
\label{sec:xndelta_proof}
In this section we prove the two propositions \ref{prop:xndelta1} and
\ref{prop:xndelta2}. These are central to
the proof of the two theorems concerning the behaviour of the estimators with
subsampled data. We start with a rough estimate on $x^\eps_{n+1} - x_n^\eps$ that we will
need for the proofs of the propositions.

\subsection{A Rough Estimate}
\begin{lemma}\label{lem:rough}
Let Assumptions \ref{a:1} hold and assume that $x^\eps (t)$, the solution of
\eqref{e:xeps_est}, is stationary. Then there exists a constant $C$, independent of
$\delta$ and $\epsilon$, such that
\begin{equation}\label{e:est_rough}
\E^{\mu^\eps} \left| x^\eps(s) - x^\eps_{n \delta} \right|^p \leq C \left( \delta^p +
\delta^{\frac{p}{2}} + \eps^p \right),
\end{equation}
for every $s \in (n \delta, (n +1 ) \delta]$ and every $p \geq 1$.
\end{lemma}
\proof
Using the same derivation that leads to \eqref{e:integr_parts},
but with $(n+1) \delta$ replaced by $s$, we have:
\begin{eqnarray}
x^\eps(s) - x_n^\eps & = & - \alpha \int_{n \delta}^{s} V'(x^\eps(s)) (1 +
\partial_y \phi(y^\eps(s))) \, ds + \sqrt{2 \sigma} \int_{n \delta}^{s} (1 +
\partial_y \phi(y^\eps(s))) \, d \beta(s)
\nonumber \\ &&- \eps \left( \phi(y^\eps (s)) - \phi(y^\eps(n\delta)) \right)
\nonumber \\  & =:&
I_{n,\delta}^1 + I_{n,\delta}^2 + I_{n,\delta}^3.
\label{eq:itophi}
\end{eqnarray}
We need to estimate the terms in \eqref{eq:itophi}. We start with $I^3_{n, \delta}$. By
Lemma \ref{l:Poisson} we have
$$\|\phi(y)\|_{L^\infty} \leq C. $$
Consequently
$$
\E^{\mu^\eps} |I_{n,\delta}^3|^p  \leq C \eps^p.
$$
To estimate $I_{n,\delta}^1$ we use again Lemma \ref{l:Poisson} to conclude that
\begin{equation}\label{e:phi_est}
\|1 + \partial_y \phi(y)\|_{L^\infty} \leq C.
\end{equation}
The above estimate, together with Assumptions \ref{a:1}, Corollary \ref{cor:moments} and
the stationarity of the process $x^\eps (t),$ give
\begin{eqnarray*}
\E^{\mu^\eps} |I_{n,\delta}^1|^p & \leq & C \delta^{p-1} \int_{n \delta}^{(n+1) \delta}
\E^{\mu^\eps} |V'(x^\eps(s))|^p \, ds
 \\ & \leq & C \delta^{p-1} \int_{n \delta}^{(n+1)\delta} \E^{\mu^\eps}
|x^\eps (s) |^{p} \, ds
\\ & \leq & C \delta^{p}.
\end{eqnarray*}
Estimate \cite[Eqn. 3.25, p. 163]{KSh91} on moments of stochastic integrals,
together with equation \eqref{e:phi_est}, enable us to conclude that
\begin{eqnarray*}
\E^{\mu^\eps} |I_{n,\delta}^2|^p & \leq & C \delta^{\frac{p}{2}-1} \int_{n
\delta}^{(n+1)\delta} \E^{\mu^\eps} |1 + \partial_y \phi(y^\eps (s))|^p \, ds
\\ & \leq & C \delta^{\frac{p}{2}}.
\end{eqnarray*}
We combine the above estimates to obtain \eqref{e:est_rough}. \qed
%
%%%%%%%%%%%%%%%%%%%%%%%%%%%%%%%%%%%%%%%%%%%%%%%%%%%%%%%%%%%%%%%%%%%%%%%%%%%%%%%%%%%%%%%%
%
\subsection{Proof of Proposition \ref{prop:xndelta1}}

From Theorem \cite[Sec. 1.3]{Freid85}, \cite[Thm. 3.4.6]{KSh91} we know that the
martingale
$$M(t):=\sqrt{2\sigma}\int_0^t  \left( 1 + \partial_y \phi(y^\eps_{s}) \right)ds$$
is equal in law to a time--changed Brownian motion,
$$M(t)=\widehat{\beta}  \left(2 \sigma \int_0^t \left( 1 + \partial_y
 \phi(y^\eps(s)) \right)^2 \, d s \right).
$$
Also the quadratic variation satisfies
$$\langle M \rangle_t =2 \sigma \int_0^t \left( 1 + \partial_y \phi(y^\eps(s))
\right)^2 \, d s \approx 2\Sigma t.$$
Indeed
\begin{eqnarray*}
\E^{\mu^{\eps}} \langle M \rangle_t &=& 2 \sigma \bbE^{\mu^{\eps}}
\int_0^t \left( 1 + \partial_y \phi(y^\eps(s))
\right)^2 \, d s
\\ & = & 2\Sigma t,
\end{eqnarray*}
where the last equality follows from equation \eqref{e:coeffs} for $d = 1$. Using these
observations we write
\begin{eqnarray*}
J_n  & := & \sqrt{2 \sigma} \int_{n \delta}^{(n+1) \delta} \left( 1 + \partial_y
\phi(y^\eps(s)) \right) \, d \beta(s) \\
&=&\sqrt 2\sigma \int_0^{(n+1)\delta}\left( 1 + \partial_y \phi(y^\eps(s)) \right) \, d
\beta(s) -\sqrt 2\sigma \int_0^{n\delta}\left( 1 + \partial_y
\phi(y^\eps(s)) \right) \, d \beta(s) \\
&=& \widehat{\beta}(2\Sigma(n+1)\delta)-\widehat{\beta}(2\Sigma n\delta)
+r_{n+1}-r_n\\
&=&\sqrt {2\Sigma \delta} \xi_n+r_{n+1}-r_n,
\end{eqnarray*}
where the $\xi_n$ are i.i.d unit Gaussian random variables and
$$r_n=\widehat{\beta}(\langle M\rangle_{n\delta})-\widehat{\beta}(2\Sigma n\delta).$$

To estimate this difference we follow the proof of \cite[Thm. 2.1]{HairPavl04}. We start
by employing the H\"{o}lder continuity of Brownian motion, together with H\"{o}lder
inequality, to estimate:
\begin{eqnarray*}
\E^{\mu^\eps} \left|  \widehat{\beta} (\langle M \rangle_{n\delta}) - \widehat{\beta} (
\E^{\mu^\eps} \langle M \rangle_{n\delta}) \right|^p   & \leq & \E^{\mu^\eps} \left|
\mbox{H\"{o}l}_{\beta}(\widehat{\beta})  \left(  \langle M \rangle_{n\delta} -
\E^{\mu^\eps} \langle M \rangle_{n\delta} \right)^{\beta} \right|^p
\\ & \leq &
\E^{\mu^\eps} \left| \mbox{H\"{o}l}_{\beta}(\widehat{\beta})  \right|^{p}
\left( \E^{\mu^\eps} \left|\langle M \rangle_{n\delta} - \E^{\mu^\eps}
\langle M \rangle_{n\delta}  \right|^{\beta  q} \right)^{\frac{p}{q}}
\\ & \leq &
C \left( \E^{\mu^\eps} \left| \int_0^{ n\delta} H(y^\eps (z)) \, dz \right|^{\beta  q}
\right)^{\frac{p}{q}},
\end{eqnarray*}
with $\beta \in \left(0, \frac{1}{2} \right)$. We have used the notation
$$
H(y) := 2 \sigma  \left( 1 + \partial_y \phi(y) \right)^2  - 2 \Sigma.
$$
We have also used the fact that, for every $\beta \in \left(0, \frac{1}{2} \right)$ and
every bounded time interval, the $\beta$--H\"{o}lder exponent of Brownian motion is
uniformly bounded with probability one. We have that
$$
\int_{\T} H(y) \, \mu(dy) = 0,
$$
where $\mu(dy)$ is defined in \eqref{e:gibbs_torus}. Since $n\delta \le T$,
Lemma \ref{lem:ito} applies
and we have that,  for $q$ sufficiently large and for $\eps$ sufficiently small,
\begin{eqnarray*}
\E^{\mu^\eps} \left| J_n - \sqrt{2 \Sigma \delta}\xi_n \right|^p   & \leq & C \left( \eps^{2 q
\beta} +   \eps^{q \beta }  \right)^{\frac{p}{q}} \\ & \leq & C
\left( \eps^{2 p \beta} + \eps^{p \beta }  \right) .
\end{eqnarray*}
This completes the proof of the proposition.
 \qed
%
%%%%%%%%%%%%%%%%%%%%%%%%%%%%%%%%%%%%%%%%%%%%%%%%%%%%%%%%%%%%%%%%%%%%%%%%%%%%%%%%%%%%%%%
%
\subsection{Proof of Proposition \ref{prop:xndelta2}}
We have
\begin{eqnarray*}
\E^{\mu^\eps} |R_2(\eps,\delta)|^p & = & \E^{\mu^\eps} \left| \int_{n \delta}^{(n +1)
\delta} \alpha V'(x^\eps(s)) \left( 1 +
\partial_y \phi(y^\eps(s)) \right) \, ds  -
 \delta A V'(x^\eps_{n \delta}) \right|^p  \\ &=&
\E^{\mu^\eps} \left| \int_{n \delta}^{(n +1) \delta} \alpha V'(x^\eps_{n\delta})\left( 1 +
\partial_y \phi(y^\eps(s)) \right) \, ds - A \int_{n\delta}^{(n+1) \delta} V'(x^\eps_{n\delta}) \,
ds \right.
\\ &&
+  \left. \alpha \int_{n\delta}^{(n+1) \delta} \Bigl(V'(x^\eps(s)) - V'(x^\eps_{n
\delta})\Bigr)\Bigl(1+
\partial_y \phi(y^\eps(s))\Bigr) \,ds \right|^p
\\ & \leq & C \E^{\mu^\eps} \left| V'(x^\eps_{n\delta})
\int_{n \delta}^{(n +1) \delta} \left( \alpha \left( 1 + \partial_y \phi(y^\eps(s))
\right) - A \right) \, ds \right|^p
\\ &&
+ \alpha ^p C \E^{\mu^\eps} \left| \int_{n\delta}^{(n+1) \delta} \Bigl( V'(x^\eps(s)) -
V'(x^\eps_{n \delta}) \Bigr) \Bigl(1+\partial_y \phi(y^\eps(s))\Bigr) \,ds \right|^p \\ &
=: & I^1_{\eps, \delta} + I^2_{\eps, \delta},
\end{eqnarray*}
where the constant $C$ depends only on $p$. We use the H\"{o}lder inequality, Assumptions
\ref{a:1}, Lemma \ref{lem:rough} and the uniform bound on $\partial_y \phi(y)$ to obtain,
for $\eps, \, \delta$ sufficiently small,
\begin{eqnarray*}
I_{\epsilon,\delta}^2  & \leq & C \delta^{p - 1} \int_{n\delta}^{(n+1) \delta}
\E^{\mu^\eps} \left| x^\eps(s) - x^\eps_{n \delta} \right|^{p} \, ds
\\ & \leq & C \delta^{p-1} \int_{n\delta}^{(n+1) \delta}(\delta^{\frac{p}{2}} +
\eps^{p}  ) \, ds \\ & \leq & C \left(\delta^{\frac{3p}{2}} + \delta^p \eps^{p}
   \right).
\end{eqnarray*}
Consequently
\begin{equation}\label{e:est_xn1}
 \left(\E^{\mu^\eps} |I^2_{\eps, \delta}| \right)^{1/p}  \leq   C(\delta^{3/2} +\delta
\epsilon).
\end{equation}
Consider now the function
$$
H(y):= \alpha \left( 1 + \partial_y \phi(y) \right) - A,
$$
From the definition of $A$ we get that
$$
\int_{\bbT} \Bigl( \alpha \left( 1 + \partial_y \phi(y) \right) - A \Bigr)\, \mu(dy) = 0.
$$
Hence, Lemma \ref{lem:ito} applies and we get
\begin{eqnarray*}
\E^{\mu^\eps} \left| \int_{n \delta}^{(n +1) \delta} \left( \alpha \left( 1 + \partial_y
\phi(y^\eps(s)) \right) - A \right) \, ds \right|^p
 & \leq & C \left(\eps^{2 p} + \eps^p \delta^p + \eps^p \delta^{p/2} \right).
\end{eqnarray*}
We combine the above estimate with \eqref{e:linbnd} and Corollary \ref{cor:moments} to obtain,
\begin{equation}\label{e:est_xn2}
 \left(\E^{\mu^\eps} |I^1_{\eps, \delta}|^p \right)^{1/p}  \leq   C \left(\eps^2 +  \eps
\delta^{1/2 } \right),
\end{equation}
for $\eps, \, \delta$ sufficiently small. The proof of the proposition follows from
estimates \eqref{e:est_xn1} and \eqref{e:est_xn2}. \qed
%
%
%%%%%%%%%%%%%%%%%%%%%%%%%%%%%%%%%%%%%%%%%%%%%%%%%%%%%%%%%%%%%%%%%%%%%%%%%%%%%%%%%%%%%%%%%%
%
%                        PROOF OF THM 1.2
%
%%%%%%%%%%%%%%%%%%%%%%%%%%%%%%%%%%%%%%%%%%%%%%%%%%%%%%%%%%%%%%%%%%%%%%%%%%%%%%%%%%%%%%%%%%
\section{Proof of Main Theorems}

Here we combine the results from the preceding two sections to complete the proofs of the
main theorems.

\subsection{Proof of Theorem \ref{thm:est_ddim}}
\label{sec:ddim}
We combine equations \eqref{e:a_est} and \eqref{e:xeps_V} to calculate
\begin{eqnarray*}
\widehat{A}(x^{\eps}) & = &  \frac{\int_0^T - \langle \nabla V(x^\eps(t)) , d x^\eps (t)
\rangle}{\int_0^T | \nabla V(x^\eps(t))|^2 \, dt}
\\& = &
\frac{\int_0^T  \left\langle - \nabla V( x^\eps(t)),  - \alpha \nabla V(x^{\eps}(t)) \,
dt -  \frac{1}{\epsilon}\nabla p \left( \frac{x^\eps(t)}{\eps} \right) \, dt + \sqrt{ 2
\sigma} \, d\beta(t) \right\rangle }{\int_0^T |\nabla V(x^\eps(t))|^2 \, dt}
\\& = &
\alpha +  \frac{\frac{1}{\epsilon}\int_0^T  \left\langle \nabla V(x^\eps(t)), \nabla
p(\frac{x^\eps(t)} {\eps} ) \right\rangle \, dt}{\int_0^T |\nabla V(x^\eps(t))|^2 \, dt}
- \sqrt{2 \sigma} \frac{\int_0^T  \left\langle \nabla V(x^{\eps}(t)) , d \beta(t)
\right\rangle}  {\int_0^T | \nabla V(x^\eps(t))|^2 \, dt}
\\
& = :& \alpha + I_1(T, \eps) - I_2(T, \eps).
\end{eqnarray*}
We will treat the terms $I_1(T, \eps)$ and $I_2(T,\eps)$ separately. We start with
$I_2(t, \eps)$. Since the stochastic integral
$$
M_T :=\int_0^T  \left\langle \nabla V(x^{\eps}(t)) , d \beta(t) \right\rangle
$$
is a continuous martingale which is null at $0$, the strong law of large numbers for
martingales \cite[p. 187]{yor} applies and we have that
$$
\lim_{T \rightarrow + \infty} \frac{M_T}{\langle M \rangle_T} = 0 \quad \mbox{a.s.}
$$
Consequently
\begin{equation}
\lim_{T \rightarrow + \infty}I_2(T,\eps) =  0 \quad \mbox{a.s.}
\label{e:i2_lim_d}
\end{equation}

Let us consider now the term $I_1(T, \eps)$. We use the ergodic theorem to deduce that
\begin{eqnarray*}
\lim_{T \rightarrow \infty} I_1(T, \eps) & = & \lim_{T \rightarrow \infty} \frac{
\frac{1}{\epsilon T} \int_0^T  \left\langle \nabla V(x^\eps(t)), \nabla p \left(
\frac{x^\eps(t)}{\eps} \right) \right\rangle \, ds}{\frac{1}{T}\int_0^T |\nabla
V(x^\eps(t))|^2 \, dt}
\\ & = &
\frac{\E^{\mu^\eps} \left(\left\langle \nabla V(x),  \frac{1}{\eps} \nabla p \left(
\frac{x}{\eps} \right) \right\rangle \right) } { \E^{\mu^\eps} | \nabla V(x)|^2} \quad
\mbox{a.s.}
\end{eqnarray*}
Now we use Proposition \ref{lem:xeps_meas_ddim} to compute
\begin{eqnarray*}
 \frac{\E^{\mu^\eps} \left(\left\langle \nabla V(x),  \frac{1}{\eps}
 \nabla p \left( \frac{x}{\eps}
\right) \right\rangle \right) } { \E^{\mu^\eps} | \nabla V(x)|^2} & = & \frac{\int_{\R^d}
\left\langle \nabla V(x),  \frac{1}{\eps} \nabla p \left( \frac{x}{\eps} \right)
\right\rangle \rho^\eps(x) \, dx }{  \E^{\mu^\eps} | \nabla V(x)|^2 }
\\ & = &
\frac{-\sigma \frac{1}{Z^\eps} \int_{\R^d} \left\langle \nabla V(x) e^{-\frac
{\alpha}{\sigma} V(x)}, \nabla \left( e^{-\frac{1}{\sigma} p(x/ \eps)} \right)
\right\rangle \, dx }{ \E^{\mu^\eps} | \nabla V(x)|^2 }
\\ & = &
\sigma \frac{ \E^{\mu^\eps} ( \Delta V(x) )}{ \E^{\mu^\eps} | \nabla
 V(x)|^2} - \alpha.
\end{eqnarray*}
In deriving the penultimate line we used an integration by parts. The weak convergence of
$\mu^\eps$ to $\mu$ (second part of Proposition \ref{lem:xeps_meas_ddim}), formula
\eqref{e:gibbs}, together with another integration by parts give
\begin{eqnarray*}
\lim_{\eps \rightarrow 0} \frac{ \E^{\mu^\eps} (\Delta V(x))}{\E^{\mu^\eps} (|\nabla
V(x)|^2)} & = & \frac{\E^{\mu} (\Delta V(x))}{\E^{\mu} (|\nabla V(x)|^2)}
\\ & = & \frac{\E^{\mu} (\Delta V(x))}{ -\frac{\sigma}{\alpha} \frac{1}{Z}
\int_{\R^d} \langle \nabla V(x) , \nabla (e^{-\frac{\alpha}{\sigma} V(x))} \rangle dx }
\\ & = & \frac{\alpha}{\sigma}.
\end{eqnarray*}
We combine the above calculations to conclude that
\begin{equation}
\lim_{\eps \rightarrow 0} \lim_{T \rightarrow \infty} I_1(T, \eps) = 0 \quad \mbox{a.s.}
\label{e:i1_lim_d}
\end{equation}
The proof of the convergence of the maximum likelihood estimator, eqn.
\eqref{e:a_est_lim} now follows from equations \eqref{e:i1_lim_d} and \eqref{e:i2_lim_d}.

The proof of the convergence of the estimator for the diffusion coefficient, eqn.
\eqref{e:sigma_est_lim}, follows from the definition of the quadratic variation, see e.g.
\cite{BasRao80}. \qed
\begin{remark}
An immediate corollary of the proof of the above theorem is that
$$
\lim_{T \rightarrow \infty} \widehat{A}(x^{\epsilon}) = \sigma \frac{ \E^{\mu^\eps} (\Delta
V(x))}{\E^{\mu^\eps} |\nabla V(x)|^2} \quad \mbox{a.s.}
$$
\end{remark}
%
%
%%%%%%%%%%%%%%%%%%%%%%%%%%%%%%%%%%%%%%%%%%%%%%%%%%%%%%%%%%%%%%%%%%%%%%%%%%%%%%%%%%%%%%%%%%%%
%
%                   PROOF OF PROPOSITION 1.1
%
%%%%%%%%%%%%%%%%%%%%%%%%%%%%%%%%%%%%%%%%%%%%%%%%%%%%%%%%%%%%%%%%%%%%%%%%%%%%%%%%%%%%%%%%%%
%
\subsection{Proof of Theorem \ref{thm:par_est_alpha}}\label{sec:thm_alpha}
We combine Proposition \ref{prop:xndelta2} and \eqref{e:integr_parts} to conclude that
$$
x^\eps_{n+1} - x^\eps_n = J_n - A\delta V'(x^\eps_n) + R(\eps, \delta),
$$
where $J_n$ is as defined in the proof of Proposition \ref{prop:xndelta1}
and, for $\eps, \, \delta$ sufficiently small and $\alpha \in (0,1)$,
\begin{equation}\label{e:R_est_1}
\left( \E^{\mu^\eps} |R(\eps, \delta)|^p \right)^{1/p} \leq
C \bigl(\delta^{3/2}+\epsilon \bigr).
\end{equation}
Notice that
$$\bbE^{\mu^{\eps}} |J_n|^2={\cal O}(\delta).$$
We combine this with formula \eqref{e:alpha_estim_1d} to obtain
\begin{eqnarray}
\widehat{A}_{N, \delta}(x^\eps) &=& A -  \frac{\sum_{n=0}^{N-1} V'(x^\eps_n)
J_n}{\sum_{n=1}^{N-1} |V'(x^\eps_n)|^2 \delta} - \frac{\sum_{n=0}^{N-1} V'(x^\eps_n)
R(\eps, \delta)}{\sum_{n=0}^{N-1} |V'(x^\eps_n)|^2 \delta} \nonumber \\ & :=& A - I_1 -
I_2 \label{e:a_j1_j2},
\end{eqnarray}
We need to control the terms $I_1$ and $I_2$. We start with $I_1$, which we rewrite in
the form
\begin{eqnarray*}
I_1 & = &  \eps^{\frac{\gamma
-\alpha}{2}}\frac{\frac{1}{\sqrt{(N\delta)}}\sum_{n=0}^{N-1} V'(x^\eps_n)
J_n}{\frac{1}{N}\sum_{n=0}^{N-1} |V'(x^\eps_n)|^2}.
\end{eqnarray*}
The central limit theorem for (discrete) martingales implies that
\begin{eqnarray*}
\lim_{N \rightarrow + \infty} \frac{1}{\sqrt{(N\delta)}} \sum_{n=0}^{N-1} V'(x^\eps_n)
J_n  & = & \frac{1}{\sqrt \delta}\mathcal{N} \left(0, \E^{\mu^{\eps}} \left(
|V'(x^{\eps}(0))|^2|J_0|^2 \right) \right)
\\ & = &
\frac{1}{\sqrt \delta}\mathcal{N} \left(0, c \, \delta \right) = c \, \mathcal{N}(0,1)
\quad \mbox{in law},
\end{eqnarray*}
for some $c$ uniform in $\epsilon \to 0$. In the above we have used the fact that
$\E^{\mu^\eps} |J_0|^2 = 2 \Sigma \delta$.

On the other hand, the ergodic theorem implies that
\begin{equation}\label{e:denom}
\lim_{N \rightarrow + \infty}\frac{1}{N}\sum_{n=0}^{N-1} |V'(x^\eps_n)|^2 = \E^{\mu^\eps}
|V(x)|^2, \quad \mbox{a.s.}
\end{equation}
Hence, by Slutsky's theorem, and remembering that $N = [\eps^{-\gamma}]$, we have that
\begin{equation}\label{e:j1_lim}
\lim_{\eps \rightarrow 0 } I_1 = 0 \quad \mbox{in law}.
\end{equation}
Consider now the term $I_2$. It can be written as
$$
I_2 = \frac{\eps^{\gamma - \alpha}\sum_{n=0}^{N-1}V'(x^\eps_n) R(\eps, \delta)
}{\frac{1}{N}\sum_{n=0}^{N-1} |V'(x^\eps_n)|^2}.
$$
The ergodic theorem implies that the denominator in the above expression converges a.s.
to a finite value. To study the numerator of the above expression we use estimate
\eqref{e:R_est_1}, together with H\"{o}lder inequality to estimate
\begin{eqnarray*}
\E^{\mu^\eps} \left| \eps^{\gamma - \alpha} \sum_{n = 0}^{N-1} V'(x^\eps_n) R(\eps,
\delta)\right| &\leq & \eps^{\gamma - \alpha} \sum_{n=0}^{N - 1} \left(
\E^{\mu^\eps}|V'(x^\eps_n)|^q  \right)^{1/q} \left( \E^{\mu^\eps} |R(\eps, \delta)|^p
\right)^{1/p}
\\ & \leq &
C \eps^{\gamma - \alpha} \sum_{n = 0}^{N - 1}
\Bigl(\E^{\mu^\eps}  |R(\eps, \delta)|^p \Bigr)^{1/p}
\\ & \leq &
C\bigl(\epsilon^{\alpha/2}+\epsilon^{1-\alpha}\bigr).
\end{eqnarray*}
In the above we have used Corollary \ref{cor:moments}, together with Assumptions
\ref{a:1}. The above calculation shows that numerator of $I_2$ converges to $0$ in $L^1$,
and hence in law. This, together with the a.s. convergence of the denominator
and Slutsky's theorem gives
\begin{equation}\label{e:j2_lim}
\lim_{\eps \rightarrow 0 } I_2 = 0 \quad \mbox{in law}.
\end{equation}
Combining \eqref{e:a_j1_j2}, \eqref{e:j1_lim} and \eqref{e:j2_lim} completes the
proof of the theorem.  \qed
%
%
%%%%%%%%%%%%%%%%%%%%%%%%%%%%%%%%%%%%%%%%%%%%%%%%%%%%%%%%%%%%%%%%%%%%%%%%%%%%%%%%%%%%%%%%%%
%
%                        PROOF OF THM 1.1
%
%%%%%%%%%%%%%%%%%%%%%%%%%%%%%%%%%%%%%%%%%%%%%%%%%%%%%%%%%%%%%%%%%%%%%%%%%%%%%%%%%%%%%%%%%%
%
\subsection{Proof of Theorem \ref{thm:par_est_sigma}}\label{sec:thm_sigma}
We combine Proposition \ref{prop:xndelta1} with \eqref{e:integr_parts} to write
the difference $x_{n+1}^\eps - x_n^\eps$  in the form
\begin{equation}
x_{n+1}^\eps - x_n^\eps = \sqrt{2 \Sigma \, \delta} \, \xi_n + \widehat{R}(\delta, \eps)
\label{e:xn_loc_1}
\end{equation}
in law, where, for $\eps, \, \delta$ sufficiently small,
\begin{equation}\label{e:R_est_3}
\left( \E^{\mu^\eps} |\widehat{R}(\eps, \delta)|^p \right)^{1/p}
\leq C \left(\delta + \eps^{\beta}\right).
\end{equation}
We substitute \eqref{e:xn_loc_1} into the formula for the estimator
\eqref{e:sigma_estim_1d} with $d = 1$ to obtain
\begin{eqnarray*}
\widehat{\Sigma}_{N, \delta}(x^{\epsilon}) &=& \Sigma \frac{1}{N} \sum_{n=0}^{N-1}
\xi_n^2 + \frac{1}{2 N \delta} \sum_{n=0}^{N-1} \left( \widehat{R}(\delta, \eps)
\right)^2 + \frac{1}{N \delta} \sum_{n=0}^{N-1} \sqrt{2\Sigma \delta}\xi_n
\widehat{R}(\delta, \eps)
\\ & =:&
\Sigma \frac{1}{N}\sum_{n=0}^{N-1} \xi_n^2 + I_1 + I_2.
\end{eqnarray*}
By the law of large numbers the first term tends almost surely to $\Sigma$ as $\epsilon
\to 0$ (which implies $N \to \infty.$) Thus it suffices to show that the remaining terms
tend to zero in law. We do this by showing that they tend to zero in $L^1.$

Note that
\begin{align*}
\bbE^{\mu^{\eps}}|I_1| & \le C \sum_{n=0}^{N-1} \bbE^{\mu^\eps}(\widehat{R}(\delta,\eps))^2\\
&=CN(\delta+\eps^{\beta})^2\\
&\le C(\delta+\epsilon^{2\beta}\delta^{-1})\\
&=C(\epsilon^{\alpha}+\epsilon^{2\beta-\alpha})\\
&=o(1),
\end{align*}
for $\alpha \in (0,1)$, since $\beta$ can be chosen arbitrarily close to $\frac12.$

Similarly
\begin{align*}
\bbE^{\mu^{\eps}}|I_2| & \le C \sum_{n=0}^{N-1} \delta^{\frac12}(\delta+\eps^{\beta})\\
&\le C(\delta^{\frac12}+\epsilon^{\beta}\delta^{-\frac12})\\
&=C(\epsilon^{\frac{\alpha}{2}}+\epsilon^{\beta-\frac{\alpha}{2}})\\
&=o(1),
\end{align*}
for $\alpha \in (0,1)$, since $\beta$ can be chosen arbitrarily close to $\frac12.$
This completes the proof.
\qed
%
%
%%%%%%%%%%%%%%%%%%%%%%%%%%%%%%%%%%%%%%%%%%%%%%%%%%%%%%%%%%%%%%%%%%%%%%%%%%%%%%%%%%%%%%%
%
%
\subsection{Proof of Theorem \ref{prop:drift_estim_2}}
Taking the limit $T \to \infty$ in \eqref{eq:alpha2} gives
$$
\lim_{T \to \infty}\tilde{A}(x^{\eps})=\widehat{\Sigma} \frac{\bbE^{\mu^{\eps}} (\Delta
V(x))}{\bbE^{\mu^{\eps}} |\nabla V(x)|^2}.
$$
Proposition \ref{lem:xeps_meas_ddim} implies that
$$\lim_{\eps \to 0}\widehat{\Sigma} \frac{\bbE^{\mu^{\eps}} (\Delta V(x))}
{\bbE^{\mu^{\eps}} |\nabla V(x)|^2}= \widehat{\Sigma} \frac{\bbE^{\mu} ( \Delta V(x)
)}{\bbE^{\mu} |\nabla V(x)|^2},$$
where $\E^\mu$ denotes expectation with respect to the invariant distribution $\rho(x)$
of the homogenized process, given by formula \eqref{e:gibbs}. An
integration by parts now gives that
$$
\E^{\mu} |\nabla V(x)|^2 = \frac{\sigma}{\alpha} \E^{\mu} ( \Delta V(x) ).
$$
Thus, the final result of our considerations is that
$$
\lim_{\eps \rightarrow 0} \lim_{T \rightarrow \infty} \tilde{A}(x^{\eps}) =
\frac{\widehat{\Sigma}}{\sigma} \alpha.
$$
\qed
%
%%%%%%%%%%%%%%%%%%%%%%%%%%%%%%%%%%%%%%%%%%%%%%%%%%%%%%%%%%%%%%%%%%%%%%%%%%%%%%%%
%
%                          CONCLUSIONS
%
%%%%%%%%%%%%%%%%%%%%%%%%%%%%%%%%%%%%%%%%%%%%%%%%%%%%%%%%%%%%%%%%%%%%%%%%
%
\section{Conclusions and Future Work}
\label{sec:conc}
The problem of parameter estimation for continuous time multiscale
diffusion processes is studied in this paper. Our goal is to accurately fit a
homogenized equation from data which has a multiscale character.
Our main conclusions are as follows:

\begin{itemize}

\item In order to estimate the drift and diffusion
coefficients accurately it is necessary to subsample.

\item There is an optimal subsampling rate, between the two
charateristic time-scales of the multiscale data.

\item The optimal subsampling rate may differ for different
parameters.

\item For gradient multiscale systems it is only necessary to estimate
the diffusion coefficient correctly, if one uses
the second estimator for the drift -- $\tilde{A}$, defined in equations \eqref{eq:alpha2} and
\eqref{e:alpha_estim_1d2}.

\end{itemize}

Both analysis and numerics are given to substantiate these claims.
Many open questions remain; we list those which seem important to us.

\begin{itemize}

\item Rough heuristics indicate that any subsampling
rate which is between the two characteristic time scales of the processes, namely
$\mathcal{O}(\eps^2)$ and $\mathcal{O}(1)$, should enable accurate  estimation
of the drift and diffusion coefficients. However our analysis works only
in the case where the subsampling is between
$\mathcal{O}(\eps)$ and $\mathcal{O}(1)$. Closing the gap between intuition
and what can be proved would be valuable.

\item Analyze other parameter estimation problems for multiscale
diffusions, not necessarily of gradient form. In particular study both
averaging and homogenization set-ups, as outlined in the introductory
section.

\item In this paper we have generated simulated multiscale data
by using a multiscale diffusion process. However this was done to
provide a convenient analytical framework. In applications it
is of interest to develop tools for characterizing the multiscale
structure of a given path  -- to estimate characteristic time--scales.
Related work has been done in \cite{FPSS03}. Further study would be
of interest.

\item Determine precisely the range of subsamplings which will give
accurate parameter estimates and optimize the subsampling rate for
accuracy.

\item Optimize the algorithm by combining estimates based on shifts of
the subsampled data -- so that information is not thrown away; this is
done in the context of econometrics and finance in
\cite{AitMykZha05b, AitMykZha05a}.

\item Analyze questions analogous to those raised here for
multidimensional multiscale processes.
\item Analyze questions analogous to those raised here
for hypoelliptic multiscale diffusions; in particular the case where the homogenized
equation is a fully elliptic first order Langevin equation which is derived from an
overdamped second-order Langevin equation.
\item Study whether there is any advantage in using random subsampling rates.
\item Study drift that depends non--linearly on the parameters to be
estimated:
$$
d x^\eps (t) = - \nabla V(x^\eps (t), \eps ; \alpha) dt + \sqrt{2 \sigma} d \beta (t).
$$
\item Parameter estimation for deterministic multiscale problems where the fast
process is a strongly mixing chaotic deterministic process.
\end{itemize}

{\it Acknowledgements} The authors are grateful to Ch. Sch\"{u}tte
for useful discussions concerning molecular dynamics, leading us to
formulate this problem. They also thank S. Olhede for useful discussions
and comments.
\bibliography{../bibtex_files/mybib}
\bibliographystyle{plain}
\end{document}